\documentclass[12pt]{extarticle}
\usepackage{amsmath, amsthm, amssymb, color}
\usepackage[colorlinks=true,linkcolor=blue,urlcolor=blue]{hyperref}
\usepackage{graphicx}
\usepackage{caption}
\usepackage{mathtools}
\usepackage{enumerate}
 \usepackage{comment}
\usepackage{verbatim}
\usepackage{tikz,tikz-cd,tikz-3dplot}
\usepackage{amssymb}
\usetikzlibrary{matrix}
\usetikzlibrary{arrows}
\usepackage{algorithm}
\usepackage[noend]{algpseudocode}
\usepackage{caption}
\usepackage[normalem]{ulem}
\usepackage{subcaption}
\usepackage{multicol}
\oddsidemargin \evensidemargin
\usepackage{makecell}
\usepackage{array}
\usepackage{enumitem}
\usepackage{xcolor}
\usepackage{listings}

 \newtheorem{theorem}{Theorem}[section]

\newtheorem{proposition}[theorem]{Proposition}
\newtheorem{lemma}[theorem]{Lemma}
\newtheorem{corollary}[theorem]{Corollary}

\theoremstyle{definition}
\newtheorem{definition}[theorem]{Definition}

\newtheorem{remark}[theorem]{Remark}

\newtheorem{example}[theorem]{Example}

\newcommand{\PP}{\mathbb{P}}
\newcommand{\RR}{\mathbb{R}}
\newcommand{\QQ}{\mathbb{Q}}
\newcommand{\CC}{\mathbb{C}}
\newcommand{\ZZ}{\mathbb{Z}}
\newcommand{\NN}{\mathbb{N}}

\DeclareMathOperator{\Cay}{Cay}

\newcommand{\boldA}{\mathfrak{A}}
\DeclareMathOperator{\diag}{diag}

\definecolor{juliakeyword}{RGB}{43, 108, 176}   % Steel Blue (function, for, end)
\definecolor{julianumber}{RGB}{0, 0, 205}       % Royal Blue (1, 2, etc.)
\definecolor{juliasymbol}{RGB}{46, 139, 87}     % Forest Green (:t, :v, :z)
\definecolor{juliaoperator}{RGB}{205, 112, 0}   % Amber/Orange (=, =>, *)

\lstdefinelanguage{JuliaOscar}{
  morekeywords={
    function, end, for, in, if, else, elseif, while, 
    return, module, using, import, export, const, struct, mutable
  },
  sensitive=true,
  morecomment=[l]{\#},
  morecomment=[n]{\#=}{=\#},
  morestring=[b]",
  otherkeywords={:t, :v, :z, =>, =},
}

\lstnewenvironment{oscarcode}[1][]
{
  \lstset{
    language=JuliaOscar,
    basicstyle=\ttfamily\small\bfseries, % Bold monospaced text to match image
    backgroundcolor=\color{white},
    keywordstyle=\color{juliakeyword},
    numberstyle=\color{julianumber},
    numbers=none,                        % No line numbers
    showstringspaces=false,
    breaklines=true,
    keywords=[2]{:t, :v, :z},
    keywordstyle=[2]\color{juliasymbol},
    keywords=[3]{=>, =},
    keywordstyle=[3]\color{juliaoperator},
    literate=
      *{0}{{{\color{julianumber}0}}}{1}
       {1}{{{\color{julianumber}1}}}{1}
       {2}{{{\color{julianumber}2}}}{1}
       {3}{{{\color{julianumber}3}}}{1}
       {4}{{{\color{julianumber}4}}}{1}
       {5}{{{\color{julianumber}5}}}{1}
       {6}{{{\color{julianumber}6}}}{1}
       {7}{{{\color{julianumber}7}}}{1}
       {8}{{{\color{julianumber}8}}}{1}
       {9}{{{\color{julianumber}9}}}{1},
    #1
  }
}{}

\title{\bf The Logarithmic A-Discriminant}
\author{Madeleine Goertz, Bernd Sturmfels and Simon Telen
 \bigskip \\
\begin{footnotesize} {\em Dedicated to the memory of Martin Helmer} \end{footnotesize}}

\date{}

\begin{document}
\maketitle

\begin{abstract} \noindent
We extend the GKZ theory of discriminants  to the setting of Euler integrals where both
 coefficients and exponents are variables. The integrand is the  likelihood function of a 
toric model  in algebraic statistics, whose number of critical points is the volume of the 
Newton polytope.  Our logarithmic $A$-discriminant vanishes when two critical points 
come together. This discriminant is a hypersurface, parametrized by the logarithmic Horn uniformization.
Hurwitz forms of toric varieties yield a degree formula, and a description of leading coefficients.
We emphasize applications and computations.
\end{abstract}

\section{Introduction}

{\em Euler integrals }\cite{MMT} are ubiquitous in mathematics and its applications. They take
the form
\begin{equation} \label{eq:eulerint}  \int_\Gamma f_z(t)^s t^\nu {\rm d} t , \end{equation}
where $\Gamma \subset \CC^d$, and one integrates
 a monomial $t^\nu = t_1^{\nu_1} t_2^{\nu_2} \cdots t_d^{\nu_d}$
times a power  of a polynomial
\begin{equation}
 \label{eq:polynomialf}
\qquad  f_z(t) \,\, = \,\, z_1 t^{a_1} \,+\,
 z_2 t^{a_2} \,+\,\cdots \, +\,z_n t^{a_n}, \qquad
  {\rm where} \,\,\,A = (a_1,a_2,\ldots,a_n) \in \ZZ^{d \times n}.
 \end{equation}
The integral \eqref{eq:eulerint} is a solution to the \emph{GKZ system} of partial differential equations 
\cite[Chapter 17]{Telen}. It is a function of the coefficients $z=(z_1,z_2,\ldots,z_n)$ and the
 exponents~$(s,\nu) = (s,\nu_1,\ldots,\nu_d)$.
 We view these as the coordinates of the product of projective spaces
 $\PP^{n-1} \times \PP^d$.
 
Limits of the integral can be studied via the
critical points of the integrand \cite[Section 2]{MMT}. These lie in
$\,\mathcal{U}_z = \big\{t \in (\CC^*)^d \,:\, f_z(t) \not= 0\bigr\}$. They
 are the solutions to the {\em critical~equations}
\begin{equation}
\label{eq:loglike}
{\rm d}_t{\rm log} \bigl( f_z(t)^s t^{\nu} \bigr)
 \,\, = \,\,  \frac{s}{f_z(t)} 
 \cdot \sum_{i=1}^d \frac{\partial f_z}{\partial t_i} {\rm d} t_i \,
 + \,\sum_{i=1}^d \frac{\nu_i }{t_i} {\rm d} t_i \quad = \quad 0.
  \end{equation}
  These  are known as {\em likelihood equations}
  in statistics, and as {\em scattering equations} in physics.

We assume that our readers are familiar with toric geometry, as presented in the textbook \cite{Telen}. We require in this section that the matroid defined by affine relations among the columns of $A$ is connected of rank $d+1$, and that $a_1,\ldots,a_n$ affinely span the lattice $\ZZ^d$.
%The matrix $A$ determines a rank $d+1$ matroid on $n$ elements.
%We now assume 
%that this matroid is connected and that $a_1,\ldots,a_n$ affinely span the lattice $\ZZ^d$.
Hence the Newton polytope $P = {\rm conv}(A)$ does not decompose as a join.
The monomials $t^{a_i}$ parametrize a $d$-dimensional
 {\em projective toric variety} $X_A$   in $\PP^{n-1}$ whose degree is
 the normalized volume $v = {\rm vol}(P)$. Our hypotheses ensure that  $v \geq 2$.
 The Euler characteristic of the very affine variety $\mathcal{U}_z$
 equals $(-1)^d\, v$, provided  the {\em principal $A$-determinant} $E_A$ does 
 not vanish at~$z$. 
 
  After clearing denominators,
  the $d$ equations in (\ref{eq:loglike}) describe a linear section of $X_A$.
  Hence our integrand $f_z(t)^s t^\nu$ has $v$ critical points, given by
  an algebraic function in $z,s,\nu$. A critical point $t \in {\cal U}_z$
  of \eqref{eq:loglike} is \emph{degenerate} if the Hessian determinant
  vanishes at $t$. Our object of study in this paper is the associated
  discriminant, which  arises from the following~result.

\begin{theorem} \label{thm:irredcodim1}
The closure of the set of parameters $(z,(s,\nu))$ for which $s \neq 0$ and (\ref{eq:loglike})
has a degenerate critical point in $\,{\cal U}_z$
  is an irreducible hypersurface~$\Lambda_A$  in the space~$\,\PP^{n-1} \times \PP^d$.
\end{theorem}

This theorem excludes the case $s = 0$, in which \eqref{eq:criteqs} is inconsistent unless $\nu = 0$ as well. 
We pass from $\,\PP^{n-1} \times \CC^d_{s \neq 0}\,$ to $\,\PP^{n-1} \times \PP^d\,$ in order to work with bihomogeneous objects.

If the affine column matroid of $A$ is not connected, then the locus $\Lambda_A$ is still 
a well-defined hypersurface, but it is reducible.
This is explained in Section \ref{sec3}; see Equation \eqref{eq:factordisconnected}.

We identify $\Lambda_A$ with its defining polynomial  in  $n\!+\!1\!+\!d$ variables $z,s,\nu$.
 This has relatively prime integer coefficients, and is thus unique up to sign.
We call $\Lambda_A$ the {\em logarithmic $A$-discriminant}. It extends the
 $A$-discriminant $\Delta_A$ of Gel'fand, Kapranov and Zelevinsky~\cite{GKZ}.

\begin{example}[Square] \label{ex:square}
Let $d=2$, $n=4$, and
$ A  =  \left [ \begin{smallmatrix} 
0 & \! 1 & \! 0 &\! 1 \\
0 &\! 0 & \! 1 & \! 1 \end{smallmatrix}  \right ]$, so 
$f_z(t )\,=\, z_1 + z_2 t_1 + z_3 t_2 + z_4 t_1 t_2$.
The surface $X_A \simeq \PP^1 \times \PP^1 $ has degree $v=2$ in $\PP^3$.
The $A$-discriminant is the $2 \times 2$ determinant, i.e.~$ \,\Delta_A = z_1 z_4 - z_2 z_3$.
By clearing denominators in (\ref{eq:loglike}),   we obtain  the bilinear equations
 \begin{equation}
 \label{eq:twoquadrics}
 \nu_1 \cdot f_z \,\,+ \,\,s \cdot t_1  \frac{\partial f_z}{\partial t_1} \,\,\, = \,\, \,
 \nu_2 \cdot f_z\,\,+ \,\,s \cdot t_2  \frac{\partial f_z}{\partial t_2} \,\,\, = \,\,\, 0. 
 \end{equation}
Their two solutions depend algebraically on the parameters
$z= (z_1,z_2,z_3,z_4)$ and $(s,\nu_1,\nu_2)$.
There is only one solution of multiplicity two when the following polynomial vanishes:
\begin{equation}
\label{eq:squarelambda}
\Lambda_A  \,\,= \,\,
 z_2^2 z_3^2 \cdot s^2
 \,+\, 
\Delta_A  \bigl(\, \Delta_A \cdot (\nu_1^2+\nu_2^2)
\,-\,2z_2z_3 \cdot s (\nu_1+\nu_2)
\,-\,2 (z_1 z_4+z_2 z_3) \cdot \nu_1 \nu_2\,\bigr). 
\end{equation}
This is our logarithmic $A$-discriminant.
Geometrically, $\Lambda_A$ is an irreducible hypersurface
of degree $(4,2)$ in $\PP^3 \times \PP^2$.
We view it as a logarithmic version of the
$2 \times 2$-determinant $\Delta_A$.
\end{example}

We now discuss the organization and main results in this paper.
Section \ref{sec2} introduces some notation and gives the formal definition of $\Lambda_A$. Its main result is 
the following characterization of all parameters $(z, s, \nu)$ for which the equations
 (\ref{eq:loglike}) exhibit the generic behavior.

\begin{theorem} \label{thm:ERL}
 There are exactly $v$ critical points in $\,\mathcal{U}_z$
 if and only if   $z$ and $(s \neq 0,\nu)$ satisfy
\begin{equation}
\label{eq:tripleproduct}
 E_A(z) \cdot R_A(s,-\nu) \cdot \Lambda_A(z,s,\nu) \,\,\not= \,\,0. 
 \end{equation}
Here, $E_A$ is the principal $A$-determinant, and
   $R_A$ defines the resonance arrangement in $\PP^d$.
\end{theorem}

The original motivation for the authors of \cite{GKZ} to introduce $A$-discriminants arose from their study of generalized hypergeometric functions, such as the integral \eqref{eq:eulerint}. In that setting, the factors $E_A(z)$ and $R_A(s,-\nu)$ are standard objects \cite[\S1.3 and \S2.9]{GKZEuler}. Theorem \ref{thm:ERL} motivates a detailed study of the third factor $\Lambda_A$.
This is the focus of the present article.

Theorem~\ref{thm:irredcodim1} is proved in
Section \ref{sec3}. In that section we derive
a parametric representation of the logarithmic $A$-discriminant
$\Lambda_A $. This is analogous to Kapranov's 
Horn uniformization~\cite{Kap} of the $A$-discriminant, to be reviewed in Proposition \ref{prop:Horn}.
Our new {\em logarithmic Horn uniformization} appears in Theorem  \ref{thm:logpara}.
A key ingredient is the configuration hypersurface of
 Denham, Schulze and Walther \cite{DSW}. In
 Lemma \ref{lem:sneaky} we write this as
   a Hadamard~product.

Section \ref{sec4} concerns the incidence correspondence in
$\, X_A \times \PP^{n-1} \times \PP^d \,$ which is
defined by the condition that the linear section (\ref{eq:loglike}) meets $X_A$ non-transversally. Its projection
into $\PP^{n-1} \times \PP^d$ contains $\Lambda_A$.
Theorem \ref{thm:projinc} states that equality holds when
the toric variety
 $X_A$ is smooth, and it
characterizes all other
  components when $X_A$ is singular.
  
In Section \ref{sec5} we express $\Lambda_A$
as a specialization of the {\em Hurwitz form} of $X_A$.
Theorem \ref{thm:fromHurwitz} gives
a formula for the bidegree of $\Lambda_A$.
We show this construction for
$A$ from~Example~\ref{ex:square}.

\begin{example} \label{ex:square2}
The two equations in (\ref{eq:twoquadrics})
are encoded in the rows of the $2 \times 4$ matrix
\begin{equation}
\label{eq:twobyfour}
Y \, =\, (y_{ij}) \,\,\, = \,\,\, \begin{small}
\begin{bmatrix}
\nu_1 & s & 0 \\
\nu_2 & 0 & s \end{bmatrix} \end{small}
\cdot \begin{footnotesize}
\begin{bmatrix} 1 & 1 & 1 & 1 \\
0 & 1 & 0 & 1 \\
0 & 0 & 1 & 1  \end{bmatrix} \end{footnotesize}
\cdot
 {\rm diag}(z_1,z_2,z_3,z_4).
\end{equation}
By \cite[Example 2.3]{Stu},
 the Hurwitz form of the square $A$ is
 the {\em $ 2 \times 2 \times 2$ hyperdeterminant}
$$
\begin{matrix} {\rm Hu}_A \,\,\, = \,\,\,
 y_{11}^2 y_{24}^2+y_{12}^2 y_{23}^2+y_{13}^2 y_{22}^2+y_{14}^2 y_{21}^2
\, +\,4 \,(y_{11} y_{14} y_{22} y_{23}+ y_{12} y_{13} y_{21} y_{24}) \phantom{xxxxxxxxxxx} \\ \quad
-2 \,(y_{13} y_{14} y_{21} y_{22}+ y_{12} y_{14} y_{21} y_{23} + y_{12} y_{13} y_{22} y_{23}
+ y_{11} y_{14} y_{21} y_{24}+ y_{11} y_{13} y_{22} y_{24}+ y_{11} y_{12} y_{23} y_{24}).
\end{matrix}
$$
The logarithmic $A$-discriminant $\Lambda_A$ is found by substituting
$Y $ into ${\rm Hu}_A$ and dividing by~$s^2$.
\end{example}

In Section \ref{sec6} we fix $z$ generic and we view $\Lambda_A$
as a polynomial in the $d+1$ variables~$(s,\nu)$.
The Newton polytope lies in a scaled $d$-simplex. We derive formulas for the
leading coefficients, as products in $\QQ[z]$  of certain resultants and mixed discriminants.
This rests on degenerations  as in \cite{BB, Stu}.
For instance, $\Delta_A$  is always a factor of $\Lambda_A$ at $s=0$, as seen in (\ref{eq:squarelambda}).
  
The concept of a  logarithmic discriminant was first introduced by
Kayser, Kretschmer and Telen \cite{KKT} in the context of hyperplane arrangements.
Section \ref{sec7} develops the common generalization of their setting and ours.
Therein, we consider likelihood functions of the~form
\begin{equation}
\label{eq:arrangement} \qquad
f(t)^s t^\nu \,\, = \,\, f_1(t)^{s_1} \cdots f_m(t)^{s_m} \, t_1^{\nu_1} \cdots t_d^{\nu_d} 
\qquad {\rm for} \quad m \geq 2.
\end{equation}
Now, $\boldA = (A_1,\ldots,A_m)$ is a tuple
of matrices.  The polynomial
$f_i(t)$ has support $A_i \in \ZZ^{d \times n_i}$
with unknown coefficients $z_{ij}$.
The case of hyperplane arrangements is
$A_1 = \cdots = A_m = ({\bf 0},{\bf e}_1,\ldots,{\bf e}_d)$. But, in contrast to
\cite{KKT}, the coefficients of the $f_i(t)$
are variables, so our logarithmic $\boldA$-discriminant
 $\Lambda_\boldA$ defines a hypersurface in
$\,\PP^{n_1-1}\times\cdots\times\PP^{n_m-1} \times \PP^{m+d-1}$.

Section \ref{sec8} features algorithms 
for computing the logarithmic $A$-discriminant $\Lambda_A$.
Our implementations reflect the state of the art in elimination theory,
which connects symbolic computation (e.g.~in {\tt Oscar} \cite{oscarbook}) 
with monodromy methods in numerical algebraic geometry.

In Section \ref{sec9} we turn to applications of the integrand
$f_z(t)^s t^\nu$ and its critical points.
In algebraic statistics, this is the likelihood function
of a scaled toric model \cite{ABB, TW}. The coefficients
$z_1,\ldots,z_n$ are scaling parameters, and the
exponents $s,\nu_1,\ldots,\nu_d$
are {\em sufficient statistics} for the data.
The maximum likelihood estimator is a function of both
of these. 

Our $\Lambda_A$ generalizes the {\em data-discriminant} \cite{RT} which 
 characterizes when the number of real critical points can change.
Another motivation comes from Feynman integrals in physics \cite{Weinzierl}.
In {\em Landau analysis} \cite{CCG, FMT, HPT, MT}, one seeks
  parameters $(z,s,\nu)$ for which (\ref{eq:eulerint}) is singular. 
However,  (\ref{eq:loglike}) shows  the scattering equations,   not the Landau equations.
When $\Lambda_A$ vanishes, the representation of all
 master integrals by a basis of Lefschetz thimbles collapses; see~\cite{ACM}.

\section{Critical Equations}
\label{sec2}

We work towards a toric formulation of our system (\ref{eq:loglike}) of $d$ rational equations in $d$ unknowns: 
\begin{equation} \label{eq:criteqs} \quad \frac{\nu_1}{t_1} + s \, \frac{\partial_1 f_z(t)}{f_z(t)} \, = \, \cdots \, = \, \frac{\nu_d}{t_d} + s \, \frac{\partial_d f_z(t)}{f_z(t)} \, = \, 0,
\quad {\rm where}\, \,\,\partial_j = \frac{\partial}{\partial t_j} .\end{equation}
For fixed $z, s, \nu$, these are constraints on points $t$ in the
very affine variety $\mathcal{U}_z = (\mathbb{C}^*)^d \setminus V(f_z)$. When $s = 0$, there are no solutions. We exclude this degenerate case and assume that $s \neq 0$. 

The exponent matrix $A$  fixes a toric compactification of ${\cal U}_z$: the closure of the image of 
\begin{equation} \label{eq:PhiA} \Phi_A: (\mathbb{C}^*)^d \, \rightarrow \, \mathbb{P}^{n-1}, \quad  t \, \longmapsto \, (t^{a_1} : t^{a_2} : \cdots : t^{a_n}) \end{equation}
is the projective toric variety $X_A \subset \mathbb{P}^{n-1}$. Our assumption that the columns of $A$ affinely generate $\mathbb{Z}^d$ ensures that $\Phi_A$ is injective. Using homogeneous coordinates $x = (x_1: \cdots: x_n)$ on $\mathbb{P}^{n-1}$, 
the monomial map $\Phi_A$ identifies ${\cal U}_z$ with the complement in $X_A$ of $n + 1$ hyperplanes: 
\[ {\cal U}_z \, \simeq \, \Phi_A({\cal U}_z) \, = \, \{ x \in X_A \, : \, x_1 x_2 \cdots x_n \, (z \cdot x) \, \neq \, 0 \}. \]
Here we write $\,z \cdot x = z_1 \, x_1 + \cdots + z_n \, x_n$. This observation is useful because, in $x$-coordinates, the equations \eqref{eq:criteqs} are linear. To make this precise, we introduce some shorthand notation.

Let $\hat{A}$ be the $(d+1) \times n$ matrix obtained by placing a row of ones on top of $A$, and set 
\begin{equation} \label{eq:notation} \hat{A}_z \, \coloneq \, \hat{A} \,{\rm diag}(z_1, \ldots, z_n) \, \, \text{ for any } z \in \mathbb{C}^n, \quad \mu \coloneq (s, - \nu). \end{equation}
Throughout the paper, we assume that the column matroid of $\hat{A}$ is connected, but this is not needed in this section.
The map $\Phi_A$ from \eqref{eq:PhiA} identifies the set of solutions of \eqref{eq:criteqs} with 
\begin{equation} \label{eq:critcirc7} {\rm Crit}^\circ \, = \, {\rm Crit}^\circ(z, s, \nu) \, = \,  \{ x \in X_A \, : \, {\rm rank} \begin{bmatrix} \hat{A}_z \, x & \mu \end{bmatrix} \, \leq \, 1 \, \,  \text{and} \, \,  \, x_1 x_2 \cdots x_n \, (z \cdot x) \, \neq \, 0 \}.\end{equation}
Using \eqref{eq:criteqs}, the rank condition ${\rm rank} \begin{bmatrix} \hat{A}_z \, x & \mu \end{bmatrix} \, \leq \, 1$ 
is written in $t$-coordinates as follows:
\begin{equation}
\label{eq:tworowscrit}
{\rm rank}  \begin{bmatrix} 
f_z(t) & t_1 \cdot \partial_1 f_z(t) &  t_2 \cdot \partial_2 f_z(t) & \cdots & t_d \cdot \partial_d f_z(t) \\
-s & \nu_1 & \nu_2 & \cdots & \nu_d  \end{bmatrix} \, \leq \, 1.
\end{equation}
We view ${\rm Crit}^\circ$ as a subscheme of $X_A$. For a geometric interpretation, it is convenient to drop the inequation
 in \eqref{eq:critcirc7}. We do this in two steps.
  First, we allow solutions on the toric boundary. That is, we study critical points on the affine open subset $z \cdot x \neq 0$ of $X_A$. We set
\[ {\rm Crit}^{\rm aff} \, = \,  {\rm Crit}^{\rm aff}(z, s, \nu) \, = \,  \{ x \in X_A \, : \, {\rm rank} \begin{bmatrix} \hat{A}_z \, x & \mu \end{bmatrix} \, \leq \, 1 \, \,  \text{and} \, \,  \,  z \cdot x \, \neq \, 0 \}.\]
Second, we obtain linear equations on $X_A$ by removing the inequation $z \cdot x \neq 0$: 
\begin{equation} \label{eq:Crit} {\rm Crit} \, = \, {\rm Crit}(z, s, \nu) \, = \, \{ x \in X_A \, : \, {\rm rank} \begin{bmatrix} \hat{A}_z \, x & \mu \end{bmatrix} \, \leq \, 1 \}.\end{equation}

It is useful to note that, since $s \neq 0$, our rank condition is equivalent to
\begin{equation} \label{eq:coeffmtx} T_{s, \nu} \, \hat{A}_z \, x = 0, \quad \text{where} \quad   T_{s, \nu} \, = \, \small \begin{bmatrix}
\nu_1 & s & 0 & \cdots & 0 \\ 
\nu_2 & 0 & s & \cdots & 0 \\ 
\vdots & \vdots & \vdots & \cdots &  \vdots \\ 
\nu_d & 0 & 0 & \cdots & s
\end{bmatrix}. \normalsize \end{equation}
At this point, we have introduced three critical loci which form a
chain of inclusions:
\begin{equation}
\label{eq:3chain}
{\rm Crit}^\circ \,\,\subseteq \,\,{\rm Crit}^{\rm aff} \,\,\subseteq \,\,{\rm Crit}.
\end{equation}
Each of these loci depends on the choice of parameters $(z,s,\nu)$.
When does equality hold?

 The factors in (\ref{eq:tripleproduct}) answer that question.
Recall from \cite{GKZ} that the {\em $A$-resultant}  ${\rm Res}_A$ is the Chow form of $X_A$.
The \emph{principal $A$-determinant} is the following specialization thereof:
 \begin{equation}
 \label{eq:EA}
 E_A \,\,=\,\, {\rm Res}_A\bigl(\,f_z,\, t_1  \partial_1 f_z, \,\ldots, \, t_d \, \partial_d f_z \bigr).
 \end{equation}
Geometrically, $E_A(z) = 0$ defines the hypersurface in $\mathbb{P}^{n-1}$ obtained by projecting the incidence $\{(x, z) \in X_A \times \mathbb{P}^{n-1} \, : \, \hat{A}_z \, x \, = \, 0\}$ to $\mathbb{P}^{n-1}$ \cite[Section 8.2]{Telen}. 
This has many irreducible factors. Each of these is the discriminant associated to a face of the polytope $P = {\rm conv}(A)$.

\begin{lemma} \label{lem:compareYX}
The right inclusion in (\ref{eq:3chain})
is strict if and only if $E_A(z) = 0$. 
\end{lemma}

\begin{proof}
If the inclusion is strict, then $T_{s,\nu} \hat{A}_z \,  x = 0$ has a solution $x $ in $X_A$ with  $z \cdot x = 0$. 
That solution  satisfies $\hat{A}_z x = 0$, which implies  $E_A(z) = 0$, by definition. Conversely, if $E_A(z) = 0$, 
then there exists a point $x \in X_A$ such that $\hat{A}_z x = 0$. In particular,
this point $x$ satisfies $z \cdot x = 0$ and $T_{s,\nu} \hat{A}_z x = 0$, and hence it
 lies in ${\rm Crit} \backslash {\rm Crit}^{\rm aff} $.  So, the inclusion is strict.
\end{proof}

Consider the convex polyhedral cone in $\RR^{d+1}$ that is spanned by the columns of $\hat A$.
Each of its facets spans a hyperplane, defined by a linear form in $\mu  = (s,-\nu)$ with
relatively prime $\ZZ$-coefficients.
Let $R_A$ be the product of these linear forms times $s$.
Thus $R_A = 0$ defines the arrangement in $\PP^d$ given by the facet hyperplanes of
the polytope $P = {\rm conv}(A)$.
Because of its role for $A$-hypergeometric functions \cite{GKZEuler}, we call this the \emph{resonance arrangement} of $A$.

The resonance arrangement appears in the next lemma, whose proof uses some more notation. For a vector $z \in \mathbb{C}^{n}$ and a face $F$ of $P$, $z_F \in \mathbb{C}^{|S_F|}$ is the restriction to $S_F = \{ i \, : \, a_i \in F \}$, and $\hat{A}_F$ is the submatrix of $\hat{A}$ whose columns are indexed by $S_F$. We write $O_F \subset X_A$ for the torus orbit of $F$. That is, $x \in O_F$ means that $x_i \neq 0$ if and only if $i \in S_F$. 

\begin{lemma}\label{lem:compareYU}
If the  left inclusion in (\ref{eq:3chain})
is strict, then $R_A(\mu) =0$,
i.e.~$(s,-\nu)$ is in the resonance arrangement.
Conversely, if $E_A(z) \neq 0$ and $R_A(\mu) = 0$, then
the inclusion is~strict.
\end{lemma}

\begin{proof}
Suppose  ${\rm Crit}^\circ \subsetneq {\rm Crit}^{\rm aff}$. 
Then ${\rm rank} \begin{bmatrix} \hat{A}_z \, x & \mu \end{bmatrix}  \leq 1$ has a solution $x$ on the toric boundary of $X_A$. 
Let $F$ be the proper face of $P$ whose torus orbit $O_F$ contains $x$.
Then  ${\rm rank} \begin{bmatrix} (\hat{A}_F)_{z_F}x_F & \mu \end{bmatrix} \leq 1$. This means that $\mu$ lies on the column span of $\hat{A}_F$, so $R_A(\mu) = 0$. 

Conversely, suppose $R_A(\mu) = 0$ and let $F$ be the minimal face of $P$, with respect to inclusion, such that $\mu$ is
in the column span of $\hat{A}_F$. We claim that there exists a solution $x$ of the equations ${\rm rank} \begin{bmatrix} \hat{A}_z \, x  & \mu \end{bmatrix} \leq 1$
in the orbit $O_F$, and that solution satisfies $z \cdot x  = z_F \cdot x_F \neq 0$. Hence, $x \in {\rm Crit}^{\rm aff} \backslash {\rm Crit}^\circ$. On $O_F$, the equations reduce to ${\rm rank} \begin{bmatrix} (\hat{A}_F)_{z_F} x_F & \mu \end{bmatrix} \leq 1$. The key insight is that these are the equations \eqref{eq:criteqs} for the toric variety $X_{A_F}$ with coefficients $z_F$. Indeed, assuming $\dim F = r$, we may assume that $F$ lies in the coordinate subspace where the last $d-r$ coordinates are zero. Then the last $d-r$ rows of $\begin{bmatrix} (\hat{A}_F)_{z_F} x_F & \mu \end{bmatrix}$ are zero, and imposing rank $\leq 1$ gives precisely the critical equations of $(\sum_{i\in S_F} z_i t^{a_i})^s t_1^{\nu_1} \cdots t_r^{\nu_r}$. There is at least one solution on $X_{A_F}$, since a general codimension-$r$ linear subspace must intersect the $r$-dimensional toric variety $X_{A_F}$. 
By Lemma \ref{lem:compareYX} and the hypothesis $E_A(z) \not= 0$, which implies $E_{A_F}(z_F) \neq 0$, there are no solutions satisfying $z_F \cdot x_F = 0$. Since $F$ is the smallest face with $\mu$ in the column span of $\hat{A}_F$, there is no solution $x_F \in X_{A_F}\setminus O_F$ which lies in a smaller torus orbit. 
Hence, ${\rm rank} \begin{bmatrix} \hat{A}_z \, x & \mu \end{bmatrix}\leq 1$ has a solution $x \in O_F \setminus \{ z \cdot x = 0 \}$. 
\end{proof}

\begin{corollary} \label{cor:ER}
Both inclusions in (\ref{eq:3chain})
are equalities if and only if $E_A(z)R_A(s, -\nu) \neq 0$. 
\end{corollary}

Corollary \ref{cor:ER} brings us close to proving Theorem \ref{thm:ERL}. Before completing the proof,
 let us formally define the main object of interest, the logarithmic $A$-discriminant. We set
\[ {\cal U} \, = \, \{ (t, z) \in (\mathbb{C}^*)^d \times \mathbb{P}^{n-1} \, : \, f_z (t) \neq 0 \}. \]
This is a family of very affine varieties ${\cal U}_z$ over $\mathbb{P}^{n-1}$. 
With $\Phi_A$ as in \eqref{eq:PhiA}, we define
 the map 
\begin{equation} \label{eq:Psi}  \Psi\, :\,  {\cal U} \, \longrightarrow \, \mathbb{P}^{n-1} \times (\mathbb{P}^d \setminus \{s = 0 \}), \quad (t,z) \, \longmapsto \, (z, \hat{A}_z \, x), \qquad
\hbox{where $x = \Phi_A(t)$.}  \end{equation}
The coordinates on $\mathbb{P}^d$ are $\mu = (s, - \nu)$. The fiber of $\Psi$ over $(z, (s, -\nu))$ is $\{z \} \times {\rm Crit}^\circ(z,s,\nu)$. In particular, $\Psi$ is generically finite of degree $v$. The \emph{ramification divisor} ${\cal R} \subset {\cal U}$ of $\Psi$ consists of the points in ${\cal U}$ at which the differential ${\rm d} \Psi$ drops rank. The defining equation of ${\cal R}$ will be given in Section \ref{sec3}. The \emph{branch locus} of $\Psi$ is the image of ${\cal R}$ under $\Psi$. This consists of all points $(s, -\nu)$ for which ${\rm Crit}^\circ(z,s, \nu)$ is either non-reduced or positive dimensional. 

\begin{definition} \label{def:logAdisc}
The \emph{logarithmic $A$-discriminant variety} is the closure in $\mathbb{P}^{n-1} \times \mathbb{P}^d$ of the branch locus of $\Psi$. Its defining equation $\Lambda_A \in \mathbb{Z}[z, s, \nu]$ is the \emph{logarithmic $A$-discriminant}.
\end{definition}

Theorem \ref{thm:irredcodim1} states that the logarithmic $A$-discriminant variety has codimension $1$
and is irreducible. We will prove this in Section \ref{sec3} by further analyzing the ramification locus ${\cal R}$. 
For now, we assume that it holds, so that the
   logarithmic $A$-discriminant $\Lambda_A$ is well-defined. We will write $\Lambda_A$ for both the logarithmic $A$-discriminant hypersurface and its defining equation. 
 
\begin{proof}[Proof of Theorem \ref{thm:ERL}]
Suppose $\# ({\rm Crit}^\circ) = v$. Then both inclusions in
(\ref{eq:3chain}) are equalities, and hence $E_A(z)R_A(s, -\nu) \neq 0$,
by Corollary \ref{cor:ER}. This tells us that 
${\rm Crit}^\circ$ is a finite scheme of length $v$.
The hypothesis $\# ({\rm Crit}^\circ) = v$ means that it is reduced,
and hence $\Lambda_A(z,s,\nu) \not= 0$. 
Conversely, suppose  (\ref{eq:tripleproduct}) holds.
Then both equalities hold in (\ref{eq:3chain}), and hence
${\rm Crit}^\circ$ is a finite scheme of length $v$.
The hypothesis $\Lambda_A(z,s,\nu) \not= 0$ implies that
$\# ({\rm Crit}^\circ) = v$.
\end{proof}

Theorem \ref{thm:ERL} strengthens
Huh's Theorem   on ML degrees 
 for scaled toric models.
The formulation in \cite{Huh} applies to the case when $z$ is fixed:
for generic $s,\nu$, the
number of critical points of $f_z(t)^s t^\nu$ in
$\mathcal{U}_z$ equals the signed Euler characteristic
$(-1)^d \chi( \mathcal{U}_z)$, and all critical points are non-degenerate.  An earlier result by
Khovanskii states that $(-1)^d \chi( \mathcal{U}_z) = v$.
Our new result characterizes
precisely what it means for $s,\nu$ to be generic in  Huh's Theorem.

To us, Theorem~\ref{thm:ERL} is a strong motivation for
studying the logarithmic $A$-discriminant $\Lambda_A$. It is
the unique irreducible factor in (\ref{eq:tripleproduct}) which depends
on both sides of the product $\PP^{n-1} \times \PP^d$.
While the other factors are well-studied in the literature, 
$\Lambda_A$ is more mysterious. 

\begin{remark} \label{rmk:dualdefective}
The logarithmic $A$-discriminant $\Lambda_A$ is always a non-constant
polynomial, even when the toric variety $X_A$ is {\em dual defective}.
This means that the $A$-discriminantal variety $\nabla_A$ is not
a hypersurface, and 
the $A$-discriminant $\Delta_A$ is $1$.
This happens for many  Segre varieties $X_A$,
including varieties of matrices of rank $1$
when the matrix format is not square.
\end{remark}

\begin{example}[The logarithmic $2 \times 3$ determinant]  \label{ex:23independence}
Let $X_A$ be the Segre threefold $\PP^1 \!\times\! \PP^2$ in $\PP^5$.
This is cut out by the $2 \!\times\! 2$ minors of a $2 \!\times \!3$ matrix
\begin{footnotesize} $\begin{bmatrix} x_1 & \! x_2 \! & x_3 \\ x_4 & \! x_5 \! & x_6 \end{bmatrix}$.
\end{footnotesize}
Here $d = 3$, $n=6$,
$$ 
A \,\, = \,\, \begin{small}
\begin{bmatrix}
0 & 1 & 0  & 0 & 1 & 0  \\
0 & 0 & 1 & 0 & 0 & 1 \\
0 & 0 & 0 & 1 & 1 & 1 \\
\end{bmatrix} , \end{small} \quad {\rm and} \,\,\quad  f_z(t) \, = \,
 z_1 \,+\, z_2 t_1 \,+\, z_3 t_2 \,+\,
 z_4 t_3 \,+\, z_5 t_1 t_3 \,+\, z_6 t_2 t_3.
$$
The triangular prism $P = {\rm conv}(A)$ has normalized volume $v=3$.
The toric variety $X_A$~is defective: its dual has codimension $2$ and hence $\Delta_A = 1$.
The principal $A$-determinant equals
\begin{equation}
\label{eq:EAex}
 E_A(z) \,\, = \,\, z_1z_2z_3z_4z_5z_6 (z_1z_5-z_2z_4) (z_1 z_6 - z_3 z_4)( z_2 z_6 - z_3 z_5). 
 \end{equation}
This has degree $12 = v (d+1)$, as expected from (\ref{eq:EA}). The resonance arrangement equals
\begin{equation}
\label{eq:RAex}
 R_A(s, -\nu) \,\, = \,\, \nu_1 \, \nu_2 \, \nu_3 \, (\nu_3 + s) \, (\nu_1 + \nu_2 + s) \, s. 
\end{equation}
All critical points of $  f(t_1,t_2,t_3)^s t_1^{\nu_1} t_2^{\nu_2} t_3^{\nu_3}$
are in $\mathcal{U}_z$ provided (\ref{eq:EAex}) and (\ref{eq:RAex}) are both non-zero.
Its coordinates can be  written in radicals in $z,s,\nu$, thanks to Cardano's formula
  for solving the cubic. There are $3$ distinct critical points
  when  the logarithmic $A$-discriminant is non-zero. We find that
  $\Lambda_A$ has $441$ terms, each 
  of degree $12$ in   $z$    and degree $4$ in $(s,\nu)$.
  In the logarithmic setting of this paper, we view
    $\Lambda_A$  as a kind of determinant for  $2 \times 3$ matrices.
  
Our theory defines such a {\em logarithmic determinant} for arbitrary rectangular matrices.
This is $\Lambda_A$ when $A$ is a product of two simplices and
$f_z(t)$ is a general bilinear form. 
Its bidegree is given in Corollary \ref{cor:logdet}.
Likewise, we can study {\em logarithmic hyperdeterminants}.
\end{example}

\section{Logarithmic Horn Uniformization}
\label{sec3}

In this section we present a parametric representation of the
logarithmic $A$-discriminant $\Lambda_A$ in $\PP^{n-1} \times \PP^d$.
This will reveal that $\Lambda_A$ is an irreducible hypersurface, as asserted in
Theorem~\ref{thm:irredcodim1}.
 Our point of departure is Kapranov's
   Horn uniformization  \cite{Kap} 
   % and \cite[Theorem 7.3.12]{Telen})
for the $A$-discriminant $\Delta_A \subset \PP^{n-1}$. We write $L$ for the row space of the
matrix $\hat A$, and $L^\perp = {\rm kernel}(\hat A)$ for its orthogonal complement in $\RR^n$.
Hence $L$ and $L^\perp$ are subspaces $\RR^n$ of dimensions
$d+1$ and $n\!-\!d\!-\!1$ respectively. We view them as subspaces
of dimensions $d$ and $n\!-\!d\!-\!2$ in $\PP^{n-1}$. We write $x \star v = (x_1 v_1 : x_2 v_2: \cdots : x_n v_n)$ 
for the Hadamard product of $x$ and $v$ in $\mathbb{P}^{n-1}$.
\begin{proposition}[Horn Uniformization] \label{prop:Horn}
The A-discriminant equals the Hadamard product
\begin{equation}
\label{eq:asin} \Delta_A \, = \, X_A \star L^\perp \,\, = \,\,
\bigl\{ x \star v \in \PP^{n-1} \,: \, x \in X_A \,\,\,{\rm and} \, \,\,v \in L^\perp \bigr\}. 
\end{equation}
This equation holds for the $A$-discriminantal variety even when this is not a hypersurface.
\end{proposition}

We now come to our extension for the logarithmic setting. Fix the reciprocal linear space
$$ L^{-1} \,\, = \,\, \overline{\bigl\{ \,v^{-1} = (v_1^{-1}: v_2^{-1}: \cdots : v_n^{-1} ) \,\, : \,\,
v \in L \cap (\RR^*)^n \,\bigr\}} \,\,\subset \,\, \PP^{n-1}. $$
Since ${\bf 1} \in L$, this is a rational  variety
of dimension $d$ in $\PP^{n-1}$.  Our main result is as follows:

\begin{theorem} \label{thm:logpara}
The logarithmic $A$-discriminant variety
is birational to $X_A \times (L^\perp \star L^{-1})$ via
\begin{equation}
\label{eq:asin2}
 X_A \times (L^\perp \star L^{-1})  \, \dashrightarrow\, \Lambda_A \, \subset \, \PP^{n-1} \times \PP^d\,,\,\,\,
(x,\phi)\, \mapsto \,\bigl(x \star \phi , \hat A\phi\bigr). 
\end{equation}
This means that  $\,s = \sum_{i=1}^n \phi_i$, $\,\nu = - A\phi$, and
$\,z= x \star \phi\,$ is a Hadamard product as in~(\ref{eq:asin}).
The image of (\ref{eq:asin2}) is a hypersurface. The inverse map
$\Lambda_A \dashrightarrow X_A \times (L^\perp \star L^{-1})$
 is given by
\begin{equation}
\label{eq:loghorninverse} x_j \,=\, 
\left(\frac{\partial \Lambda_A}{\partial z_j} \right)^{\!\! -1} \!\! \!\cdot
\biggl( \sum_{i=1}^d a_{ij} \frac{\partial \Lambda_A }{\partial \nu_i}
\,-\,\frac{\partial \Lambda_A}{\partial s}  \biggr) \quad
{\rm for}\,\,\, j = 1,\ldots,n
\,\quad {\rm and} \,\quad
\phi \,=\, x^{-1} \star z.
\end{equation}
\end{theorem}

The partial derivatives in (\ref{eq:loghorninverse}) are polynomials in $z,s,\nu$.
These define a rational map that recovers $x_j$ whenever $(z,(s,-\nu))$ 
lies on the image hypersurface $\Lambda_A$.
Our inversion formula (\ref{eq:loghorninverse}) is a variant of Theorem 1.5 in Section 1.1.B on page 16 in \cite{GKZ},
which is a manifestation of the Biduality Theorem for projective varieties. In the toric context, this
maps a general point on the $A$-discriminant $\Lambda_A$ to the
associated point on the toric variety $X_A$, i.e.~it inverts the Horn uniformization,
provided $X_A$ is not dual defective. See \cite[Proposition 7.1.11]{Telen}.

We note that Theorem \ref{thm:irredcodim1} is an immediate corollary of Theorem \ref{thm:logpara}. 
Before giving the proof of Theorem~\ref{thm:logpara},
here is an example that shows the inversion formula (\ref{eq:loghorninverse}) in action.

\begin{example}[Square]
Fix $d=2$, $n=4$, and revisit Examples \ref{ex:square} and \ref{ex:square2}.
We use the coordinates $b$ on  $L^\perp \simeq \RR \,(1,-1,-1,1)$ and
 $(c_1,c_2,c_3)$ on  $L \simeq \RR^3$. Then
the map (\ref{eq:asin2}) equals
\begin{equation}
\label{eq:invfun1} \begin{matrix} z_1 \,=\, t_0 \cdot \frac{b}{c_1}\,,\,\,\,
 z_2 \,=\, -t_0 t_1 \cdot \frac{b}{c_1 + c_2}\,,\,\,\,
 z_3 \,=\, -t_0 t_2 \cdot \frac{b}{c_1 + c_3}\,,\,\,\,
 z_4 \,=\, t_0 t_1 t_2 \cdot \frac{b}{c_1 + c_2 + c_3}\,, \medskip \\
 s = \frac{b}{c_1}-\frac{b}{c_1 + c_2} -\frac{b}{c_1 + c_3} + \frac{b}{c_1 + c_2 + c_3},\,\,
 \nu_1 = \frac{b}{c_1 + c_2} - \frac{b}{c_1 + c_2 + c_3},\,\,
 \nu_2 = \frac{b}{c_1 + c_3} - \frac{b}{c_1 +  c_2 +  c_3} . \end{matrix}
 \end{equation}
  The image is an irreducible hypersurface of degree $(4,2)$ in $\PP^3\! \times\! \PP^2$.
 The defining polynomial $F = \Lambda_A$ of this hypersurface is displayed in~(\ref{eq:squarelambda}).
 We write $ F_s , F_{\nu_1},F_{\nu_2}, F_{z_1},F_{z_2},F_{z_3},F_{z_4}$ for the partial derivatives of $F$.
 These are polynomials in $z,s,\nu$.  The inversion formula  (\ref{eq:loghorninverse})~reads
\begin{equation}
\label{eq:invfun2} 
 x_1 \,=\, -\frac{F_s}{ F_{z_1}} , \,\, \,
 x_2 \,=\, \frac{ F_{\nu_1} - F_s}{F_{z_2}} ,\,\,\,
x_3 \,=\, \frac{ F_{\nu_2} - F_s}{F_{z_3}} ,\,\,\,
x_4 \,=\, \frac{F_{\nu_1} + F_{\nu_2}  - F_s}{ F_{z_4}} .
\end{equation}
 The right hand sides in  (\ref{eq:invfun2}) are complicated. But, after
  substituting (\ref{eq:invfun1}),  they simplify to
 $$ x_1 \,=\, t_0\,,\quad
x_2 \,=\, t_0 t_1 \,,\quad
x_3 \,=\, t_0 t_2\, ,\quad
x_4 \,=\, t_0 t_1 t_2 .$$
From this we recover $(c_1,c_2,c_3,b)$ up to scale, by the Hadamard product $\phi \,=\, x^{-1} \star z$.
 \end{example}
 
The proof of Theorem \ref{thm:logpara} requires a better understanding of the rational variety $L^\perp \star L^{-1}$. Our connectedness assumption on the affine matroid of $A$ ensures that this is the hypersurface $\{ \phi \in \mathbb{P}^{n-1} :  \det \hat{A}_\phi \hat{A}^T = 0 \}$. This was called the \emph{configuration hypersurface} of $\hat{A}$~in~\cite{DSW}. 
 
\begin{lemma} \label{lem:sneaky}
The Hadamard product $L^\perp \star L^{-1}$ is a rational hypersurface 
of degree $d+1$ in $\PP^{n-1}$.
Its defining polynomial in $\phi$ is the determinant of the 
$(d+1) \times (d+1)$ matrix $ {\hat A}_\phi {\hat A}^T$.
\end{lemma}

\begin{proof}
The Main Theorem in \cite{DSW} says that 
the configuration hypersurface ${\cal X} = \{ \phi \in \PP^{n-1}:
\det {\hat A}_\phi {\hat A}^T = 0 \}$
 is reduced and irreducible. Here we use that
  the matroid of $\hat{A}$ is connected.
  Let $u \in \mathbb{C}^n$ be in the right kernel of $\hat{A}$, and $v = \hat{A}^T w \in (\mathbb{C}^*)^n$ for some $w \in \mathbb{C}^{d+1} \backslash \{0\}$. 
Consider the point $\phi = u \star v^{-1}$ in $L^\perp \star L^{-1}$. 
By construction, the row vector $w$ lies in the right kernel of
$ {\hat A}_\phi {\hat A}^T = \hat A {\rm diag}(\phi) {\hat A}^T $.
Hence the determinant of this matrix is zero.
Thus, the image of the map $L^\perp \times L^{-1} \dashrightarrow \mathbb{P}^{n-1}$ 
given by $(u, v^{-1}) \mapsto u \star v^{-1}$ is contained in ${\cal X}$. 

We must show that the image of $L^\perp \times L^{-1} \dashrightarrow \mathbb{P}^{n-1}$ is a dense subset of ${\cal X}$. By \cite[Main Theorem]{DSW}, the smooth points of ${\cal X}$ are 
${\cal X}_{\rm sm} = \{ \phi \in \mathbb{P}^{n-1} \, : \, {\rm rank}( \hat{A}_\phi \hat{A}^T) = d \}$. 
We define
\[ {\cal X}^\circ \, = \, \bigl\{\, \phi \in {\cal X}_{\rm sm} \, : \, \hat{A}^T w \in (\mathbb{C}^*)^n \text{ for each } w \in \ker( \hat{A}_\phi\hat{A}^T )\setminus \{ 0 \} \,\bigr\}.\]
We claim that ${\cal X}^\circ$ is dense in ${\cal X}$. This follows from the fact that the incidence variety 
${\cal Y} = \{(\phi, w) \in \mathbb{P}^{n-1} \times \mathbb{P}^d \, : \, \hat{A}_\phi \hat{A}^T \, w = 0 \}$ 
is birational to ${\cal X}$ via the $\phi$-coordinate projection \cite[Proposition 3.16]{BDSW}. 
The $w$-coordinate projection ${\cal Y} \rightarrow \mathbb{P}^d$ is surjective.
Also, $\hat{A}^T w \in (\mathbb{C}^*)^n$ for generic $w$, because the matroid of $\hat{A}$ has no loops.
This implies that ${\cal X}^\circ$ is dense in ${\cal X}$. 

If $\phi \in {\cal X}^\circ$, then ${\hat A}_\phi {\hat A}^T$ has rank $d$.
Let $w$ be a vector that spans its right kernel. Since $\phi \in {\cal X}^\circ$, this satisfies $\hat{A}^T w \in (\mathbb{C}^*)^n$. Then
$\phi $ is given by $u \star (w^T \hat A)^{-1}$ for some $u \in \ker \hat{A}$. Hence, our map reaches the dense subset ${\cal X}^\circ$, and each $\phi \in {\cal X}^\circ$ admits an inverse.
\end{proof}

\begin{proof}[Proof of Theorem \ref{thm:logpara}]
By Definition \ref{def:logAdisc},
$\Lambda_A = \overline{\Psi({\cal R})}$, where ${\cal R}$ is the ramification divisor of the map $\Psi$ in \eqref{eq:Psi}. Throughout the proof, let $x = \Phi_A(t)$ be the image of $t$ under the monomial map in \eqref{eq:PhiA}. A calculation shows that  ${\cal R} = \{ (t, z) \in {\cal U} \, : \, \det \hat{A}_{z \star x} \hat{A}^T \, = \, 0 \}$. This connects our theory  with the configuration hypersurface
of Denham, Schulze and Walther~\cite{DSW}.

We claim that ${\cal R}$ is irreducible of dimension $d+n-2$. For this, notice that 
\[ {\cal R} \, \subseteq \, {\cal R}' \, = \,  \{(t,z) \in  (\mathbb{C}^*)^d  \times \mathbb{P}^{n-1} \, : \, \det \hat{A}_{z \star x} \hat{A}^T \, = \, 0 \}. \]
The righthand side is an irreducible hypersurface. Indeed, after applying the coordinate change $( t,z) \mapsto (t, z \star x)$, ${\cal R}'$ is given by $\{ (t, \phi) \in (\mathbb{C}^*)^d \times \mathbb{P}^{n-1}\, : \, \det \hat{A}_\phi \hat{A}^T = 0\}$, which is irreducible by Lemma \ref{lem:sneaky}. Since both ${\cal R}$ and ${\cal R}'$ have dimension $d + n-2$, we have
%, because $\det \hat{A}_{z \star x} \hat{A}^T = 0$ does not imply that $f_z(t) = z_1x_1 + \cdots + z_n x_n = 0$. In $(t,\phi)$ coordinates one sees this from the fact that $\{\det \hat{A}_\phi \hat{A}^T = 0\}$ is not contained in $\sum_i \phi_i = 0$. 
%Therefore, we have 
${\cal R}' = \overline{\cal R}$. 

Since $\Lambda_A = \overline{\Psi({\cal R})}$, we have now shown that $\Lambda_A$ is irreducible of codimension at least one. Moreover, by Lemma \ref{lem:sneaky}, ${\cal R}$ is parametrized by $(t,z) = (t,x^{-1} \star \phi)$, where $\phi \in L^\perp \star L^{-1}$. Composing this with $\Psi$ gives the logarithmic Horn uniformization \eqref{eq:asin2}. The formula in \eqref{eq:asin2} equivalently uses $x$ instead of $x^{-1}$. This does not change the image of the map. 

To show that $\Lambda_A$ is a hypersurface, we prove that the restriction of $\Psi$ to a dense open subset of ${\cal R}$ has finite fibers. We extend $\Psi$ to a rational map $\Psi_{\rm rat}:  (\mathbb{C}^*)^d \times \mathbb{P}^{n-1} \dashrightarrow \mathbb{P}^{n-1} \times \mathbb{P}^d$. 
Let ${\cal V} \subset \mathbb{P}^{n-1} \times \mathbb{P}^d$ be the open subset given by $s \, E_A(z) R_A(s,-\nu) \neq 0$. We claim that $\Psi_{\rm rat}^{-1}({\cal V}) \cap {\cal R}'$ is dense in ${\cal R}'$. Let us first exclude the possibility that $E_A(z)$ vanishes identically on $\Psi_{\rm rat}({\cal R}')$. If this were the case, then the image of the coordinate projection ${\cal R}' \rightarrow \mathbb{P}^{n-1}$ would be contained in $\nabla_{A_F}$ for some face $F$ of ${\rm conv}(A)$. However, we claim that ${\cal R}' \rightarrow \mathbb{P}^{n-1}$ is dominant. Indeed, since $v \geq 2$, for generic $z \in \mathbb{P}^{n-1}$, 
the determinant of $\hat{A}_{z \star x} \hat{A}^T$ is a Laurent polynomial in $t$ with at least two terms, 
so the fiber of ${\cal R}' \rightarrow \mathbb{P}^{n-1}$ over $z$ is non-empty. 

Second, we must show that $\Psi_{\rm rat}({\cal R}')$ cannot be contained in the hyperplane arrangement given by $R_A(s,-\nu) = 0$. This would imply that $\{\det \hat{A}_\phi \hat{A}^T = 0 \}$ is contained in a hyperplane, which contradicts Lemma \ref{lem:sneaky}. We conclude that $\Psi^{-1}({\cal V}) \cap {\cal R}$ is dense in ${\cal R}$. We already know 
from Corollary \ref{cor:ER} that the map
$ \Psi^{-1}({\cal V}) \cap {\cal R} \rightarrow {\cal V}$  given by $\Psi$ has finite fibers.
Therefore, we 
conclude that $\Psi_{|{\cal R}}: {\cal R} \rightarrow \Lambda_A$ is generically finite, and $\dim \Lambda_A = \dim {\cal R} = d + n - 2$. 

It remains to prove the inversion formula \eqref{eq:loghorninverse}. 
We compute a tangent vector $\tau$ to the logarithmic $A$-discriminant by differentiating
the parametrization \eqref{eq:asin2}.
We obtain
\[ \tau \, = \, (x \star \dot{\phi} + \dot{x} \star \phi, \, \hat{A} \, \dot{\phi}). \]
This expresses $\tau$ in terms of a tangent vector $\dot{\phi}$ to the hypersurface $L^\perp \star L^{-1}$ and a tangent vector $\dot{x}$ to the toric variety $X_A$. 
For a generic point $(z, (s, -\nu)) \in \Lambda_A$, let $(x, \phi)$ be a point in the pre-image of \eqref{eq:asin2} and let $u \in \mathbb{C}^{d+1}$ be a generator of the one-dimensional kernel of $\hat{A}_\phi \hat{A}^T$. We show that the normal vector to $\Lambda_A$ at $(z, (s, -\nu)) \in \Lambda_A$ is given by $(\hat{A}_{x^{-1}}^T u, -u)$:
\[ (\hat{A}_{x^{-1}}^T u, -u) \cdot \tau \, = \, u^T \hat{A}_{x^{-1}}(x \star \dot{\phi} + \dot{x} \star \phi) - u^T \hat{A} \, \dot{\phi} \, = \, u^T \hat{A}_{x^{-1}}(\dot{x} \star \phi) \, = \, 0. \]
The last equality uses the familiar fact $\dot{x} \in {\rm Row}( \hat{A}_x )$ for a tangent vector to $X_A$. We have now shown that $(\hat{A}_{x^{-1}}^T u, -u)$ is proportional to the gradient vector of $\Lambda_A$ at $(z, (s, -\nu))$, so
\[ u \, = \, c \cdot \Big ( \frac{\partial \Lambda_A}{\partial s}, -\frac{\partial \Lambda_A}{\partial \nu_1}, \ldots, -\frac{\partial \Lambda_A}{\partial \nu_d} \Big ) \quad \text{and} \quad - c \cdot \frac{\partial \Lambda_A}{\partial z_i} \, = \, x_i^{-1} \, (\hat{a}_i^T u).  \] 
This completes the proof of Theorem \ref{thm:logpara}
and hence of Theorem \ref{thm:irredcodim1}.
\end{proof}

We next discuss what happens when the matroid of $A$ is not connected. We examine the case 
where $\hat{A}$ has~two connected components, from which the general case can be extrapolated. 
There exists a matrix $Q$ such that  $Q \, \hat{A}$ is block~diagonal after permuting columns:
\begin{equation}
\label{eq:directsum}
 Q \, \hat{A} \,\, = \,\, \begin{bmatrix}
\hat{A}^{(1)} & {\bf 0} \\ 
{\bf 0} & \hat{A}^{(2)}
\end{bmatrix}.
\end{equation}
The first row of $\hat{A}^{(i)} \in \mathbb{Z}^{(d^{(i)}+1)\times n^{(i)}}$ is the all-ones vector. We use the notation
\[ Q \, 
\mu
 \, = \, \begin{bmatrix}
\mu^{(1)}\\ \mu^{(2)}
\end{bmatrix} , \quad \text{where} \quad \mu \, = \, \begin{bmatrix}
s \\ -\nu
\end{bmatrix} \, \text{ and } \, \mu^{(i)} \, = \, \begin{bmatrix}
s^{(i)} \\ -\nu^{(i)}
\end{bmatrix}. \] 
The parameters $(s^{(i)}, -\nu^{(i)}) \in \mathbb{P}^{d^{(i)}}$ are linear combinations of the original exponents $s, \nu$. We claim that, assuming $s^{(i)} \neq 0$ for $i = 1, 2$, the critical equations decouple into those for the individual toric varieties $X_{A^{(i)}}$. This is directly seen from the formulation ${\rm rank} \begin{bmatrix}
\hat{A}_z \, x & \mu 
\end{bmatrix} = {\rm rank} \begin{bmatrix}
Q \hat{A}_z \, x & Q \mu 
\end{bmatrix} \leq 1$. There is a one-to-one correspondence between \eqref{eq:criteqs} and the product
\begin{equation} \label{eq:split} {\rm Crit}^\circ_{A^{(1)}}(z^{(1)}, s^{(1)}, \nu^{(1)}) 
\times {\rm Crit}^\circ_{A^{(2)}}(z^{(2)}, s^{(2)}, \nu^{(2)}).\end{equation}
Here $z^{(i)}$ denotes the subvector of $z$ indexed by the columns of $A^{(i)}$.  

These observations indicate that the space of parameters is now naturally a product 
\begin{equation} \label{eq:splitparameters} \mathbb{P}^{n^{(1)}-1}_{z^{(1)}} \times \mathbb{P}^{d^{(1)}}_{\mu^{(1)}} \times \mathbb{P}^{n^{(2)}-1}_{z^{(2)}}  \times \mathbb{P}^{d^{(2)}}_{\mu^{(2)}}.\end{equation}
The scheme ${\rm Crit}^\circ_{A^{(i)}}$ lives in ${\cal U}_{z^{(i)}}$, a $d^{(i)}$-dimensional very affine open subset of $X_{A^{(i)}}$. Define
\[ \Psi \, : \, \prod_{i = 1, 2} \mathbb{P}^{n^{(i)}-1}_{z^{(i)}} \times {\cal U}_{z^{(i)}} \, \longrightarrow \, \prod_{i = 1, 2} \mathbb{P}^{n^{(i)}-1}_{z^{(i)}} \times (\mathbb{P}^{d^{(i)}}_{\mu^{(i)}} \setminus \{ s^{(i)} = 0 \})   \]
analogously to \eqref{eq:Psi} on each factor. The closure of the branch locus of $\Psi$ in \eqref{eq:splitparameters} is the logarithmic $A$-discriminant $\Lambda_A$. By definition, it is a hypersurface whose equation satisfies 
\begin{equation} \label{eq:factordisconnected} \Lambda_A(z, s, \nu) \, = \, \Lambda_{A^{(1)}}(z^{(1)}, s^{(1)}, \nu^{(1)}) \cdot \Lambda_{A^{(2)}}(z^{(2)}, s^{(2)}, \nu^{(2)}). \end{equation}
 Here we use the convention that $\Lambda_{A^{(i)}} = 1$ if $X_{A^{(i)}}$ has degree one, or if $d^{(i)} = 0$. 

A different choice for $Q$ only changes $\Lambda_A$ by a coordinate change on the factors~$\mathbb{P}^{d^{(i)}}_{\mu^{(i)}}$. This ambiguity can be remedied by substituting $\mu = \hat{A} \, u$ in $\Lambda_A$. This uses $u = (u_1, \ldots, u_n)$ instead of the $d+1$ exponents $(s, \nu)$, since $\Lambda_A(z, u)$ only depends on $X_A$; see
Remark~\ref{rmk:dependsonly}.

\begin{example}[Join of two conics] \label{ex:jointwoconics}
Let $d=3,n=6$ and consider the matrices
\begin{equation}
\label{eq:tetrahedron2}
 \hat{A} \, = \,  \begin{small} \begin{bmatrix}
1 & 1 & 1 & 1 & 1 & 1 \\
0 & 1 & 2 & 0 & 0 & 0 \\
0 & 0 & 0 & 1 & 1 & 1 \\
0 & 0 & 0 & 0 & 1 & 2 \\
 \end{bmatrix}, \quad Q \, = \, \begin{bmatrix}
1 & 0 & -1 & 0 \\ 0 & 1 & 0 & 0 \\ 0 & 0 & 1 & 0 \\ 0 & 0 & 0 & 1
\end{bmatrix}. \end{small}
\end{equation}
The matrix $Q$ reveals the connected components of $A$, and $\Lambda_A$ has two irreducible~factors: 
\begin{equation}
\label{eq:twohypersurfaces}
\begin{matrix}
 8z_1z_3s\nu_1 + 4z_1z_3\nu_1^2 + 8z_1z_3\nu_1\nu_2 - z_2^2s^2 - 2z_2^2s\nu_1 - 2z_2^2s\nu_2 - z_2^2\nu_1^2 - 2z_2^2\nu_1\nu_2 - z_2^2\nu_2^2,  \smallskip \\ {\rm and} \qquad
  8z_4z_6\nu_2\nu_3 - 4z_4z_6\nu_3^2 + z_5^2\nu_2^2 - 2z_5^2\nu_2\nu_3 + z_5^2\nu_3^2.
\end{matrix}
\end{equation}
These two factors are 
$\Lambda_{A^{(1)}}$ and 
$\Lambda_{A^{(2)}} $, but after a linear change of coordinates.
The first hypersurface is  parametrized by 
the logarithmic Horn map for $\hat{A}^{(1)} \, = \, \left [ \begin{smallmatrix} 1 & 1 & 1 \\ 0 & 1 & 2 \end{smallmatrix} \right ]$,
modified to
\[ z_1 \, = \, t_0 \cdot \frac{b}{c_1}, \quad z_2 \, = \, t_0t_1 \cdot \frac{-2b}{c_1+c_2}, \quad z_3 \, = \, t_0t_1^2 \cdot \frac{b}{c_1+2c_2}, \quad \begin{bmatrix} s^{(1)} \\  -\nu^{(1)} \end{bmatrix} \, = \, \begin{bmatrix}
1 & 1 & 1\\ 0 & 1 & 2
\end{bmatrix} \begin{bmatrix}
\frac{b}{c_1}\smallskip \\ \frac{-2b}{c_1+c_2 }\smallskip \\ \frac{b}{c_1+2c_2}
\end{bmatrix}\!,\]
where $s^{(1)} = s - (-\nu_2)$ and $-\nu^{(1)} = -\nu_1$. These relations between $(s, -\nu)$ and $(s^{(1)}, -\nu^{(1)})$ are given by the 
rows of $Q$.  It is instructive to examine 
the map (\ref{eq:asin2})  for the matrix $\hat{A}$ in~(\ref{eq:tetrahedron2}).
\end{example}

\begin{remark} If $\hat{A}$ is not connected then
 the image of the logarithmic Horn uniformization  is irreducible of codimension $\geq 2$.
 In Example \ref{ex:jointwoconics} it is
     the zero set of the two factors of $\Lambda_A$.
\end{remark}

We end with a characterization of dual defectivity which will be useful in later sections. 

\begin{lemma} \label{lem:dualdefect}
The toric variety $X_A$ is dual defective if and only if the determinant $\det A_v A^T$ vanishes for each $v \in \ker (\hat{A})$. That is, the configuration hypersurface of $A$ contains $\ker( \hat{A})$. 
\end{lemma}

\begin{proof} 
The dimension of the $A$-discriminantal variety $\nabla_A$ is seen from the rank of the Jacobian matrix of its parametrization,
which is the Horn uniformization (Proposition \ref{prop:Horn}). 
Let $B$ be an $(n-d-1)\times n$ matrix whose rows form a basis for $L^\perp$. Here $x_i=\Phi_A(t)$ and $v_i=(B^Tc)_i$ is the $i$-th entry of $B^Tc$. Then the Horn uniformization is given by 
$$(t_1,\dots,t_d,c_1,\dots,c_{n-d-1}) \, \longmapsto \,  \bigl( t^{a_1}(B^Tc)_1,\dots,t^{ a_n}(B^Tc)_n \bigr).$$
This parametrizes the affine cone over $\nabla_A$ in $\mathbb{C}^n$.
The $n \times (n-1)$ Jacobian matrix of this parametrization has rank $\dim \nabla_A + 1$. This equals the rank of
$\diag(x)M$ where $M=\begin{bmatrix}
A_v^T & B^T
\end{bmatrix}.$ This matrix was obtained by computing toric derivatives with respect to the toric variables $t_j$.
Since $x\in O_P$, $\diag(x)$ is invertible, so the Jacobian has the same rank as $M$, and $X_A$ is dual defective if and only if ${\rm rank}(M) < n-1$. This happens if and only if the row span of $A_v$ contains a nonzero vector in the row span of $B$, which means there exists a nonzero vector $w \in \mathbb{C}^d$ satisfying $\hat{A} A_v^T w = 0$. Since the first row of $\hat{A} A_v^T$ is zero, this is equivalent to ${\rm rank}(\hat{A} A_v^T) = {\rm rank} (A_v A^T) < d$, which concludes the proof of the lemma. 
\end{proof}

\section{Projecting the Incidence Correspondence}
\label{sec4}

In \eqref{eq:Crit}, the likelihood equations are given by a linear section of the $d$-dimensional projective toric variety $X_A$.
The codimension-$d$ linear space depends on $z, s, \nu$, and we write this as
\[ L_{z,s,\nu} \, = \, \{ x \in \mathbb{P}^{n-1} \, : \, T_{s,\nu} \hat{A}_z \, x \, = \, 0 \}. \]
This motivates us to study the following incidence~correspondence:
\[ {\cal I} \, = \, \Big \{ (x, z, \nu) \in X_A \times \mathbb{P}^{n-1} \times \mathbb{C}^d   \, : \, \begin{matrix} \text{$x$ is a non-isolated or non-reduced}\\ \text{  point of the intersection } L_{z, 1, \nu} \cap X_A  \end{matrix} \Big \}. \]
Here we set $s = 1$ to avoid the degenerate case $s = 0$. Consider the coordinate projection $\pi_1: {\cal I} \rightarrow \mathbb{P}^{n-1} \times \mathbb{C}^d$. The closure of $\pi_1({\cal I})$ in $\mathbb{P}^{n-1} \times \mathbb{P}^{d}$ contains $\Lambda_A$ by definition, but there could be other components. Our goal is to describe all components of this projection.

\begin{theorem} \label{thm:projinc}
If $X_A$ is regular in codimension one, meaning $\dim {\rm Sing}(X_A) \leq d-2$, then the hypersurface $\pi_1({\cal I})$ has only the following irreducible components: 
\begin{itemize}
\item[(a)] the logarithmic $A$-discriminant hypersurface $\Lambda_{A}$, 
\item[(b)] the $A_F$-discriminant hypersurface $\nabla_{A_F} \times \mathbb{C}^d$ 
 for each face $F $ of $P= {\rm conv}(A)$ such that $X_A$ is singular along the orbit $O_F$ and $X_{A_F}$ is 
 not dual defective (cf.~Remark \ref{rmk:dualdefective}).
\end{itemize}
In particular, if  the toric variety $X_A$ is smooth, then $\pi_1({\cal I}) = \Lambda_A$.
\end{theorem}

\begin{proof}
The locus $\pi_1({\cal I})$ is defined by a single polynomial, obtained by specializing the \emph{Hurwitz form} of $X_A$ to the matrix $T_{s, \nu} \hat{A}_z$ \cite{Stu}. This is the topic of Section \ref{sec5}. Here, we use this observation to conclude that $\pi_1({\cal I})$ is an equidimensional hypersurface. By Theorem \ref{thm:ERL}, 
its irreducible components are among the components of
$\{(z,\nu) : E_A(z) \cdot R_A(1,-\nu) \cdot \Lambda_A(z,1,\nu) = 0\}$. 
We split the rest of the proof into four parts, to be presented in Propositions \ref{prop:IP1}-\ref{prop:IF1}.
\end{proof}

Our strategy is based on a decomposition of ${\cal I}$ according to the strata of $X_A$, and according to whether $z \cdot x$ vanishes or not. For each face $F$ of the polytope $P$,  we define 
\[ {\cal I}_{F,0} \, = \, \{ (x, z, \nu) \in {\cal I} \, : \, x \in O_F \text{ and } z \cdot x = 0 \}, 
\,\,\, {\cal I}_{F,1} \, = \, \{ (x,z, \nu) \in {\cal I} \, : \, x \in O_F \text{ and } z \cdot x \neq 0 \}. \]
 As before,  $O_F$ is the torus orbit of $X_A$ corresponding to $F$, which has dimension $r = \dim F$.  The propositions below do not require that $X_A$ is regular in codimension one, but they imply Theorem \ref{thm:projinc} under that assumption. The following statement holds by definition.

\begin{proposition} \label{prop:IP1}
The closure of $\pi_1({\cal I}_{P,1})$ in $\mathbb{P}^{n-1} \times \mathbb{P}^{d}$ equals $\Lambda_A$. 
\end{proposition}

We now investigate the projections of the other strata. 

\begin{comment}
\begin{proposition} \label{prop:IP0}
The projection $\pi_1({\cal I}_{P,0})$ has codimension at least two in $\mathbb{P}^{n-1} \times \mathbb{P}^{d}$. 
\end{proposition}

\begin{proof}
Notice that $z \cdot x = 0$ implies $\hat{A}_z x = 0$ on ${\cal I}$. For $x \in O_P$, the tangent space $T_xX_A$ is the projectivization of the row span of $\hat{A}_x$. For $x$ to be a non-isolated or non-reduced solution, the tangent space should intersect the linear space $A_z x = 0$ in a line. Affinely, this means that the kernel of $\hat{A}_{z \star x} \hat{A}^T$ has dimension two. 
The equations for ${\cal I}_{P,0}$ are
\[ x \in O_P, \quad \hat{A}_z \, x \,= \, 0 \quad \text{and} \quad {\rm corank} (\hat{A}_{z \star x} \hat{A}^T) \geq 2. \]
From $\hat{A}_z x = 0$ we see that the image of $\pi_1({\cal I}_{P,0})$ in $\PP^{n-1}$ is contained in $\nabla_A$. This implies 
our assertion when $X_A$ is dual defective. Modulo $\hat{A}_z x = 0$, the corank condition can equivalently be expressed as $\det( A_{z \star x} A^T) = 0$. That determinant vanishes identically for $z \star x \in \ker( \hat{A})$ if and only if $X_A$ is dual defective, by Lemma \ref{lem:dualdefect}. Hence, if $X_A$ is not dual defective, then ${\cal I}_{P,0}$ has dimension $d + n-3$, which is an upper bound for the dimension of its projection. 
\end{proof}
\end{comment}

\begin{proposition} \label{prop:IFsing}
Suppose that $X_A$ is singular along the torus orbit $O_F$. We have 
\[ \overline{\pi_1 ({\cal I}_{F,0})} \, = \,  \nabla_{A_F} \times \mathbb{C}^d\quad \text{and} \quad \overline{\pi_1({\cal I}_{F,1})} \, = \,  \mathbb{P}^{n-1} \times -{\rm aff}_\mathbb{C}(F) .\]
\end{proposition}

\begin{proof}
The equations for ${\cal I}_{F,0}$ are given by $\hat{A}_z \, x = 0$ and $ x \in O_F$. Indeed, since $X_A$ is singular along $O_F$, any isolated solution to these equations is non-reduced. These are precisely the incidence equations for the discriminant $\nabla_{A_F}$, which proves the first claim. 

The equations for ${\cal I}_{F,1}$ are $x \in O_F$ and $x \in L_{z,1,\nu}$. Again, an isolated solution to these equations is automatically non-reduced. The condition $x \in L_{z,1,\nu}$ is equivalent to the vector $\hat{A}_z x$ being a multiple of $(1,-\nu)$.
Since $x \in O_F$, this means that $-\nu \in {\rm aff}_{\mathbb{C}}(F)$. Conversely, if $-\nu \in {\rm aff}_{\mathbb{C}}(F)$, then a solution $x \in O_F \cap L_{z,1,\nu}$ exists for generic $z$. We have seen this in the proof of Lemma \ref{lem:compareYU}. This establishes the second claim: $\overline{\pi_1({\cal I}_{F,1})}  =   \mathbb{P}^{n-1} \times -{\rm aff}_\mathbb{C}(F)$.
\end{proof}

\begin{remark} 
Proposition \ref{prop:IFsing} implies that, if $X_A$ is singular along an orbit $O_F$ of codimension one, then the divisor $\mathbb{P}^{n-1} \times -{\rm aff}_\mathbb{C}(F)$ belongs to $\pi_1({\cal I})$. This is not listed in Theorem \ref{thm:projinc}. 
\end{remark}

\begin{proposition} \label{prop:IF0}
Suppose that $X_A$ is smooth along the torus orbit $O_F$. Then
the projection $\pi_1({\cal I}_{F,0})$ has codimension $\geq 2$ in $\mathbb{P}^{n-1} \times \mathbb{C}^{d}$. 
\end{proposition}

\begin{proof}
The equations that define ${\cal I}_{F,0}$ are as follows:
\begin{equation} \label{eq:IF0} x \in O_F, \quad \hat{A}_z \, x \, = \, 0, \quad \text{and}
 \quad \dim \bigl( \mathbb{P} \ker (T_{1,\nu} \hat{A}_z) \cap T_x X_A \bigr) \,\geq\, 1. \end{equation}
Here $T_x X_A$ is the $d$-dimensional tangent space of $X_A$ at $x$. If ${\cal T} \in \mathbb{C}^{n \times (d+1)}$ is a matrix whose columns span $T_x X_A$, then the last condition is equivalent to $\,{\rm rank}(T_{1,\nu}\hat{A}_z \, {\cal T}) \leq d-1$. 

Let ${\cal J} \subset O_F \times \mathbb{P}^{n-1} \times \mathbb{C}^d$ be the irreducible $(n+d-2)$-dimensional variety defined by $(\hat{A}_F)_{z_F} \, x_F = 0$. Note that, if $X_{A_F}$ is dual defective, then $\overline{\pi_1({\cal J})} = \nabla_{A_F}\times\CC^d$ has codimension $\geq 2$ and the proposition holds. We may therefore assume that $X_{A_F}$ is not dual~defective.

We will show that the locus of points $(x,z,\nu) \in O_F\times \mathbb{P}^{n-1} \times \mathbb{C}^d$ satisfying \eqref{eq:IF0} has dimension at most $n+d-3$, so that its projection to $\mathbb{P}^{n-1}\times\CC^d$ has codimension at least two. 
The first step is to determine the matrix ${\cal T}$. Let $r = \dim F$ and suppose, without loss of generality, that ${\rm aff}_{\mathbb{C}}(F) =  \{\, q \in \mathbb{C}^d \, : \, q_{r+1} =  \cdots = q_{d} = 0\, \}$. After permuting columns, 
\begin{equation} \label{eq:blockformA} A \, = \, \begin{bmatrix}
\tilde{A}_F & \times & \times & \cdots & \times & \times \\
\bf{0} & \bf{1} & \bf{0} & \cdots & \bf{0} & \times\\
\bf{0} & \bf{0} & \bf{1} & \cdots & \bf{0} & \times \\
\vdots & \vdots & \vdots & \ddots & \vdots & \times \\
\bf{0} & \bf{0} & \bf{0} & \cdots & \bf{1} & \times
\end{bmatrix}.\end{equation}
The first block column is the matrix $A_F$ containing the columns of $A$ which lie on $F$. Its first $r$ rows are $\tilde{A}_F$. The next $d-r$ block columns consist of the columns $a_i$ of $A$ for which $a_{r+j,i} = 1$ and $a_{r+k,i} = 0$ for $k \in \{1, \ldots, d-r\} \setminus \{j\}$. For each $j \in \{1, \ldots, d-r\}$, at least one such a column $a_i$ exists by the assumption that $X_A$ is smooth along $O_F$. The first $d$ rows of these block columns are represented by the symbol $\times$, as these entries will play no role. Similarly, the last column of $\times$'s represents the remaining columns of $A$. 

Every point $x \in O_F$ can be obtained as a one-parameter limit of the following form: $x = \lim_{\epsilon \to 0} \Phi_A(t_1, \ldots, t_r, \epsilon \,  t_{r+1}, \ldots, \epsilon \,  t_d)$. In particular, $x_F = \Phi_{A_F}(t_1, \ldots, t_r)$. The tangent space $T_{\Phi_A(t)}X_A$ is spanned by the rows of $\hat{A} \, {\rm diag}(\Phi_A(t))$. Hence, the tangent space $T_x X_A$ for $x = \lim_{\epsilon \to 0} \Phi_A(t(\epsilon))$ is the limit of the row span of $\hat{A} \, {\rm diag}(\Phi_A(t(\epsilon))$ for $\epsilon \to 0$.

We calculate that limit. The first row of $\hat{A} \, {\rm diag}(\Phi_A(t(\epsilon))$ is simply $\Phi_A(t(\epsilon))$, so its limit is $x$. The next $r$ rows converge to the rank-$r$ matrix $\tilde{A}_F \, {\rm diag}(\Phi_{A_F}(t))$, with zero columns appended on the right side. Finally, the limit of the last $d-r$ rows is zero. To compute the limit of their span, we divide the $(r+j)$-th row by $\epsilon t_{r+j}$. This does not alter the row span for $\epsilon \neq 0$. Taking the limit after dividing by $\epsilon t_{r+j}$, the non-zero entries in row $r+j$ are in the position of the $\bf{1}$-vectors in \eqref{eq:blockformA}. We denote these non-zero entries by $v_{r+j}(t)$:  
\begin{equation} \label{eq:calT} T_x X_A \, = \, \mathbb{P} {\rm Row}  \begin{bmatrix} \Phi_{A_F}(t) & \bf{0} & \bf{0} & \cdots & \bf{0} & \bf{0} \\ 
\tilde{A}_F \, {\rm diag}(\Phi_{A_F}(t)) & \bf{0} & \bf{0} & \cdots & \bf{0} & \bf{0} \\
\bf{0} & v_{r+1}(t) & \bf{0} & \cdots & \bf{0} & \bf{0}\\
\bf{0} & \bf{0} & v_{r+2}(t) & \cdots & \bf{0} & \bf{0} \\
\vdots & \vdots & \vdots & \ddots & \vdots & \bf{0} \\
\bf{0} & \bf{0} & \bf{0} & \cdots & v_d(t) & \bf{0}
\end{bmatrix}.\end{equation}
The matrix in this display is ${\cal T}^T$. The matrix product $\hat{A}_z {\cal T}$ has the following block form: 
\begin{equation} \label{eq:AzT} \hat{A}_z {\cal T} \, = \, \begin{bmatrix}
z \cdot x & ( (\tilde{A}_F)_{z_F} \, x_F)^T & \times \\
(\tilde{A}_F)_{z_F}\, x_F & (\tilde{A}_F)_{z_F \star x_F} \tilde{A}_F^T & \times \\
{\bf 0} & {\bf 0} & {\rm diag}(z^{(1)} \cdot v_{r+1}(t), \ldots, z^{(d-r)} \cdot v_{d}(t)) 
\end{bmatrix}.\end{equation}
Here $z^{(j)}$ denotes the entries of $z$ corresponding to the $(j+1)$-st block column in \eqref{eq:blockformA}.

The condition ${\rm rank}(T_{1, \nu} \hat{A}_z {\cal T}) \leq d-1$ means that either ${\rm rank}(\hat{A}_z {\cal T}) \leq d-1$ or $(1, -\nu) \in {\rm im} \, \hat{A}_z {\cal T}$. We show that neither of these closed conditions holds identically on the irreducible $(n+d-2)$-dimensional variety ${\cal J}$.
This ensures that \eqref{eq:IF0}  has dimension at most $n+d-3$. 

From \eqref{eq:AzT} we see that, after imposing $(\hat{A}_F)_{z_F} \, x_F = 0$, the condition ${\rm rank}(\hat{A}_z {\cal T}) \leq d-1$~is
\begin{equation} \label{eq:corank2condition} \det ((\tilde{A}_F)_{z_F \star x_F} \tilde{A}_F^T) 
\cdot \prod_{j = 1}^{d-r} z^{(j)} \cdot v_{r+j}(t) \, = \, 0.\end{equation}
This does not hold identically on ${\cal J}$. Indeed, the fact that $\det ((\tilde{A}_F)_{z_F \star x_F} \tilde{A}_F^T)_{|{\cal J}}  \not \equiv 0$ is implied by the fact that $X_{A_F}$ is not dual defective, see Lemma \ref{lem:dualdefect}. Also $z^{(j)} \cdot v_{r+j}(t)$ does not vanish identically
because the equations $ (\hat{A}_F)_{z_F} x_F \, = \,0$ for ${\cal J}$ do not involve the variables~$z^{(j)}$. 

Finally, to see that $(1, -\nu) \in {\rm im} \, \hat{A}_z {\cal T}$ does not hold identically on ${\cal J}$, 
we note that the first column of \eqref{eq:AzT} is zero on ${\cal J}$.
Therefore, the containment $(1, -\nu) \in {\rm im} \, \hat{A}_z {\cal T}$ imposes 
at least one non-trivial linear condition on $\nu$. The equations $ (\hat{A}_F)_{z_F} x_F \, = \,0$ 
for ${\cal J}$ do not involve the variables $\nu$, so these linear conditions do not hold identically on ${\cal J}$. 
\end{proof}

\begin{proposition} \label{prop:IF1}
Suppose that $X_A$ is smooth along the torus orbit $O_F$. The
 codimension of the projection $\,\pi_1({\cal I}_{F,1})\,$ in $\,\mathbb{P}^{n-1} \times \mathbb{C}^{d}\,$
 is at least $\,{\rm codim}(F) + 1 = d-r+1$, where $\dim(F)=r$. \end{proposition}

\begin{proof}
The equations for the incidence ${\cal I}_{F,1}$ are 
$\, x \, \in \, O_F, \, T_{1, \nu} \hat{A}_z \,x \, = \, 0\,$ and $\, \det ( \hat{A}_z \, {\cal T}) \, = \, 0$.
Here ${\cal T}$ is the transpose of the matrix  in \eqref{eq:calT}. The columns of ${\cal T}$ span the $d$-dimensional tangent space to $X_A$ at the point $x$. 
The conditions $x \in O_F$ and $T_{1,\nu} \hat{A}_z \, x = 0$ imply  $\nu \in - {\rm aff}_{\mathbb{C}}(F)$. Assuming
  the form \eqref{eq:blockformA} for $A$, we have $\nu = (\tilde{\nu},\bf{0})$; the last $d-r$ coordinates are zero. This proves that $\pi_1({\cal I}_{F,1})$ has codimension at least ${\rm codim}(F) = d-r$. We  consider   \[{\cal J} \,  = \,  \{(x, z, \nu) \in  O_F\times \mathbb{P}^{n-1} \times -{\rm aff}_{\mathbb{C}}(F) \, : \, T_{1, \nu} \hat{A}_z \, x \, = \, 0 \} \quad \subset \, X_A \times \mathbb{P}^{n-1} \times \mathbb{C}^d .\]
A standard fiber argument for the coordinate projection ${\cal J} \rightarrow O_F\times -{\rm aff}_{\mathbb{C}}(F) $ shows that ${\cal J}$ is irreducible of dimension $r+n-1$. Since ${\cal I}_{F,1} = \{ (x, z, \nu) \in {\cal J} \, : \, \det ( \hat{A}_z \, {\cal T}) = 0 \}$, it suffices to show that $\det ( \hat{A}_z \, {\cal T})$ does not vanish identically on ${\cal J}$. Indeed, this implies that $\dim \pi_1({\cal I}_{F,1}) \leq \dim {\cal I}_{F,1} \leq r+n-2$. 
From \eqref{eq:AzT} we learn that, for $(x,z,\nu) \in {\cal J}$, we have 
\[ \hat{A}_z \, {\cal T} \, = \, \begin{bmatrix}
1 & -\tilde{\nu}^T & \times \\
 -\tilde{\nu}  & (\tilde{A}_F)_{z_F \star x_F} \tilde{A}_F^T & \times \\
{\bf 0} & {\bf 0} & {\rm diag}(z^{(1)} \cdot v_{r+1}(t), \ldots, z^{(d-r)} \cdot v_{d}(t)) 
\end{bmatrix}.\]
The lower right block has determinant $\prod_{j = 1}^{d-r} z^{(j)} \cdot v_{r+j}(t)$.
This  is not identically zero on ${\cal J}$ because the variables $z^{(j)}$ do not appear in the defining equations for ${\cal J}$. We also have 
\[ \det \begin{bmatrix}
1 & -\tilde{\nu}^T  \\ -\tilde{\nu}  & (\tilde{A}_F)_{z_F \star x_F} \tilde{A}_F^T  \end{bmatrix}_{|{\cal J}} \, \not \equiv \, 0.\]
Indeed, if this determinant vanishes, then $(1, \tilde{\nu}, z_F) \in \Lambda_{\tilde{A}_F}$. This concludes the proof.
\end{proof}

\section{Specializing the Hurwitz Form}
\label{sec5}

The projective toric variety $X_A$ has dimension $d$ and degree $v \geq 2$ and it lives in $\PP^{n-1}$.
We continue to assume that the singular locus of $X_A$ has codimension $\geq 2$ in $X_A$.
Under that assumption, the article \cite{Stu} defines the {\em Hurwitz form} ${\rm Hu}_A$ of $X_A$.
By the Cayley trick, ${\rm Hu}_A$
is the $A^{(d)}$-discriminant where $A^{(d)} = A \times \Delta_{d-1}$.
Here $A^{(d)}$ is a matrix with $2d$ rows and $nd$ columns.
Namely, abbreviating ${\bf 1} = (1,\ldots,1)$ and ${\bf 0} = (0,\ldots,0)$, this is the~matrix
\begin{equation}
\label{eq:A(d)}
 A^{(d)} \,\,= \,\,\begin{small} \begin{bmatrix}
\,{\bf 1} & {\bf 0} & \cdots & {\bf 0} \,\, \\
\,{\bf 0} & {\bf 1} & \cdots & {\bf 0} \,\, \\
\, \vdots & \vdots & \ddots & \vdots \,\, \\
\,{\bf 0} & {\bf 0} & \cdots & {\bf 1} \,\, \\
\,A & A & \cdots & A \,\, \\
\end{bmatrix}. \end{small}
\end{equation}
Geometrically, this is the direct product of the given polytope $P = {\rm conv}(A)$ with the simplex
$ {\rm conv} \{{\bf e}_1,{\bf e}_2,\ldots,{\bf e}_d \}$.
By \cite[Proposition 2.4]{Stu},  the Hurwitz form ${\rm Hu}_A$ is the
{\em mixed discriminant} \cite{CCDRS} for a system of $d$ Laurent polynomial equations that are supported on~$A$.
The coefficients of this polynomial system form a $d \times n$ matrix $Y = (y_{ij})$
of unknowns.  Then ${\rm Hu}_A$ is a polynomial in the $d \times d$ minors of $Y$.
These are the primal Pl\"ucker coordinates in~\cite{Stu}. The Hurwitz form vanishes when our system of equations has a singular solution in $(\mathbb{C}^*)^d$.

We define the {\em Hurwitz degree} $h = h(A)$ to be
the degree of ${\rm Hu}_A$ as a polynomial in the~$\binom{n}{d}$
Pl\"ucker coordinates. As a polynomial in the entries of $Y$,
the Hurwitz form ${\rm Hu}_A$ has degree $d h$,
and it is homogeneous of degree $h$  in each row of $Y$.
From \cite[Theorem 1.1]{Stu}, we infer:

\begin{remark} We have
$h = 2v+2g-2 $, where $g$ is the sectional genus of the toric variety  $X_A$.
This  is the  genus of the curve obtained by intersecting
$X_A$ with a general subspace of codimension $d-1$ in $\PP^{n-1}$.
If $d = 2$ then $g$ is the number of interior lattice points in the polygon $P$.
Khovanskii \cite{Kho} gives a combinatorial formula for the sectional genus~$g$ for $d \geq 3$.
\end{remark}

We now state the main result in this section.
Recall that ${\rm Sing}(X_A)$ has codimension~$\geq 2$.

\begin{theorem} \label{thm:fromHurwitz}
The logarithmic $A$-discriminant $\Lambda_A$ is a factor in the
Hurwitz form ${\rm Hu}_A$ with primal Stiefel coordinates $Y$ replaced by the entries of the matrix $T_{s,\nu} \hat{A}_z$ in \eqref{eq:coeffmtx}.
The other non-constant factors are $s^{(d-1)h}$ and powers of the discriminants of singular torus orbits of $X_A$.
The bidegree of $\Lambda_A$ is $(b,h)$ where $b \leq dh$.
Moreover, $b = dh$ if and only if $X_A$ is~smooth.
\end{theorem}

We already saw this theorem in action  in Example \ref{ex:square2}.
The proof will be given below. %The substitution $Y = T_{s, \nu} \hat{A}_z$ means that we evaluate the Hurwitz form at the system of~equations 
%\[ \nu_i \, f - s \, \partial_i f \, = \, 0, \quad, i = 1, \ldots, d. \]

\begin{example}[Tact Invariant] \label{ex:tact1}
Let $d=2,n=6$ and $A = 
\left [ \begin{smallmatrix} 
0 & 0 & 1 & 0 & 1 & 2\\
0 & 1 & 0 & 2 & 1 & 0
 \end{smallmatrix} \right ]
$. The toric variety $X_A$ is the
Veronese surface in $\PP^5$. This has degree $v=4$ 
and sectional genus $g=0$. Hence the Hurwitz degree is
$h = 6$. The Hurwitz form ${\rm Hu}_A$ is the classical {\em tact invariant}
\cite[Example 2.3]{Stu},
 a polynomial with $3210$ monomials of degree $12 = 6+6$ in  
 primal Stiefel coordinates.

We find the logarithmic $A$-discriminant  by evaluating the
tact invariant at the $2 \times 6$~matrix
\begin{equation}
\label{eq:specialY}
 Y \,\,\, = \,\,\,T_{s,\nu} \, \hat{A}_z \,\,=\,\,
\begin{bmatrix}
\nu_1 & s & 0 \\
\nu_2 & 0 & s \end{bmatrix}
\cdot \begin{bmatrix} 1 & \!\! \cdots\!\! & 1  \\ & A & \end{bmatrix} \cdot {\rm diag}(z_1,z_2,z_3,z_4,z_5,z_6).
\end{equation}
This matrix represents the two quadratic equations
$\,
\nu_1 f_z + s t_1 \frac{\partial f_z}{\partial t_1} =
\nu_2 f_z + s t_2 \frac{\partial f_z}{\partial t_2} = 0 $.
The specialized tact invariant has an extraneous factor $s^6$. This is because all three 
$2 \times 2$ minors of the $2 \times 3$ matrix on the left are divisible by $s$.
By removing the extraneous factor $s^6$, we obtain the
logarithmic $A$-discriminant $\Lambda_A$, which has $1197$ monomials  of bidegree $(12,6)$.

More generally, let $f_z(t)$ be a bivariate polynomial of degree $k$, so
$X_A$ is the Veronese surface in $\PP^{n-1}$ where $n = \binom{k+2}{2}$.
The Hurwitz degree is  $h = 2 k^2 + 2 \binom{k-1}{2} - 2 = 3k(k-1)$.
 The logarithmic $A$-discriminant $\Lambda_A$ is an irreducible
 bihomogeneous  polynomial, which has degree $6k(k-1)$ in the 
 $n$ coefficients $z$ and degree $3k(k-1)$ in the exponents $(s,\nu_1,\nu_2)$.
\end{example}

\begin{example}[Triangular Prism] \label{ex:23independenceb}
We continue Example \ref{ex:23independence}, so $d=3,n=6,v=3$.
The Hurwitz form of $A = \Delta_1 \times \Delta_2$ is
the discriminant for the $6 \times 18$ matrix $A^{(3)} =
\Delta_1 \times \Delta_2 \times \Delta_2$.
Therefore ${\rm Hu}_A$
 is the hyperdeterminant of
format $2 \times 3 \times 3$, which is the sum of
$16749$ monomials of degree $12$ in the entries
 of a $3 \times 6$ matrix $Y$.
To compute $\Lambda_A$, we evaluate ${\rm Hu}_A$ at
the special matrix $Y$  described in Theorem \ref{thm:fromHurwitz},
similar to (\ref{eq:specialY}), and we then divide by $s^8$.
\end{example}

The following example shows that $b = {\rm deg}_z(\Lambda_A)$ is less than $dh$
   when $X_A$ is singular.
   
\begin{example}
Let $d=2$, $n=5$ and $A =  \left [\begin{smallmatrix}
0 & 1 & -1 & 0 & 0 \\
0 & 0 & 0 & 1 & -1 \\
\end{smallmatrix} \right ]$. The convex hull of $A$ is a diamond with area $v=4$. It contains one interior lattice point, so $g = 1$ and $h=8$. The primal Stiefel coordinates form a $2 \times 5$ matrix $Y=(y_{ij})$.
The Hurwitz form is a polynomial of degree $16$  with $1720$ terms
in the $10$ variables $y_{ij}$. The specialization  in Theorem~\ref{thm:fromHurwitz} is here given~by
\[
Y \,\,\,=\,\, \, \begin{bmatrix}\,\,
\nu_1 z_1 & \nu_1 z_2 + s z_2 & \nu_1 z_3 - s z_3 & \nu_1 z_4 & \nu_1 z_5 \,\, \\ \,\,
\nu_2 z_1 & \nu_2 z_2 & \nu_2 z_3 & \nu_2 z_4 + s z_4 & \nu_2 z_5 - s z_5 \,\,
\end{bmatrix}. 
\]
This is obtained as $T_{s,\nu} \hat{A}_z$, similar to equation
(\ref{eq:twobyfour}). The result of this specialization is $\,{\rm Hu}_A(Y) = s^8 z_2^2 z_3^2 z_4^2 z_5^2 \cdot \Lambda_A$, where $\Lambda_A$ is an irreducible polynomial with $104$ terms. Its bidegree is $(8,8)$, which is less than the bound
$(dh,h) = (16,8)$ in Theorem~\ref{thm:fromHurwitz}. The toric surface $X_A \subset \PP^4$ has four singular points
which give rise to the extraneous factor $z_2^2z_3^2z_4^2z_5^2$.
\end{example}

\begin{proof}[Proof of Theorem~\ref{thm:fromHurwitz}]
For any $d \times n$ matrix $Y$ of rank $d$, the
Hurwitz form ${\rm Hu}_A$ vanishes at $Y$
if and only if the linear space $\{Y \cdot x = 0\}$
intersects the toric variety $X_A$ non-transversally.
We apply this to the  matrix $Y = T_{s,\nu} \, \hat{A}_z $
which represents our critical equations in \eqref{eq:coeffmtx}.
Each maximal minor of $Y$ is a polynomial of
bidegree $(d,d)$ in the unknowns $z$ and $(s,\nu)$. Hence,
the specialized Hurwitz form ${\rm Hu}_A(Y)$ 
has bidegree $(dh,dh)$ in  $z$ and $(s,\nu)$.
Moreover, each maximal minor of $T_{s,\nu}$ is
divisible by $s^{d-1}$. The same holds for the maximal minors of
$Y = T_{s,\nu} \, \hat{A}_z $. This implies that
$s^{(d-1)h}$ divides ${\rm Hu}_A(Y)$. From now,
instead of ${\rm Hu}_A(Y)$, we investigate its factor
$\, {\rm Hu}^*(Y) = s^{(1-d)h} \cdot {\rm Hu}_A(Y)$.
This has
bidegree $(dh,h)$ in $z$ and~$(s,\nu)$.

The zero set of ${\rm Hu}^*_A(Y) $  in $\PP^{n-1} \times (\PP^d \setminus \{ s = 0 \})$ is
precisely the hypersurface $\pi_1(\mathcal{I})$ in Theorem~\ref{thm:projinc}.
Hence ${\rm Hu}^*_A(Y)$ is a product of the polynomials described 
in Theorem~\ref{thm:projinc}, each appearing with a positive exponent, and with possibly an extra factor $s^M$. 
The factor (a) in ${\rm Hu}^*_A(Y)$ is a power  $(\Lambda_A)^N$ 
of the logarithmic $A$-discriminant. The factors in (b)
are the \emph{discriminants of singular orbits} in Theorem~\ref{thm:fromHurwitz}.
These are polynomials in $z$, and they do not occur when $X_A$ is smooth.
To conclude our assertion, we show that $M = 0$ and $N = 1$.

To prove this, we  specialize the coordinates of $z$ in $\Lambda_A$ to be
generic such that $E_A(z) \neq 0$. We show that the result is a polynomial of degree $h$ in $s, \nu$, which implies that $M = 0, N = 1$ since ${\rm Hu}^*_A(Y) = G(z) \cdot s^M (\Lambda_A)^N$ has bidegree $(dh, h)$. As in Section \ref{sec3}, we
consider the moment map from the toric variety $z \star X_A$
onto the projective space $\PP^d$ with coordinates $(s,\nu)$. 
This is the restriction of the linear map defined by $\hat A$. It defines a linear system $z \star X_A \dashrightarrow \mathbb{P}^d$, and $E_A(z) \neq 0$ ensures that there are no basepoints. The fiber is $z\star {\rm Crit}(z,s,\nu)$. 
Let $L \simeq \PP^1$ be a generic line in $\PP^d$. By Bertini's
Theorem, the preimage of $L$ in $z\star X_A$ is a smooth projective curve $C$.
Here, we are using that $z$ is generic and that ${\rm codim}({\rm Sing}( X_A)) \geq 2$.

The genus of the smooth curve $C$ is the sectional genus $g$ of
the toric variety $X_A$. This is because $C$ avoids the singular locus of $X_A$, so the Koszul complex resolves its structure sheaf and the Euler characteristic of $C$ equals that of a generic section of $X_A$. 
Consider the restricted moment map $C \rightarrow L$.
This map is $v$-to-$1$. 
The ramification divisor consists of $h$
points on $C$, where $h = 2v+2g-2$, by the Riemann-Hurwitz Theorem.
Each ramification index is $2$ because the configuration hypersurface is reduced.
The $h$ ramification points
project to $h$ distinct branch points on $L$. The fact that these images are distinct follows from Theorem \ref{thm:logpara}, 
which implies that, at a generic point of $\Lambda_A$, there is precisely one degenerate critical point.
The branch points are the intersection points of $L$ with the
hypersurface in $\PP^d$ that is described in
 part (a) of Theorem~\ref{thm:projinc}.
We conclude that this hypersurface is defined by
a squarefree polynomial $H$ of degree $h$ in $(s,\nu)$.
The $N$th power of $H$ factors ${\rm Hu}^*_A(Y)$,
now regarded as a polynomial in $(s,\nu)$, since $z$
had been specialized generically.
The latter polynomial has degree $h$.
Hence they are equal, and $M = 0, N=1$ as desired.
\end{proof}

We have shown that $\Lambda_A = {\rm Hu}_A^*(Y)$,
and this has bidegree $(dh,h)$, whenever $X_A$ is smooth.
As an application, we give a
formula and table for the bidegrees of logarithmic determinants.

\begin{table}[h!]
\centering
\begin{tabular}{|c||c|c|c|c|c|c|}
 \hline
  $ d_1 \,\backslash \, d_2 $ & 2 & 3 & 4 & 5 & 6 & 7 \\
 \hline\hline
  $2$    & $(4,2)$   & $ (12, 4) $ & $ (24, 6) $ & $ (40, 8) $ & $ (60, 10) $ & $ (84, 12) $ \\
    $3$    &  &  $ (48, 12) $ & $ (120, 24) $ & $ (240, 40) $ & $ (420, 60) $ & $ (672, 84) $ \\
      $4$    & & &  $ (360, 60) $ & $ (840, 120) $ & $ (1680, 210) $ & $ (3024, 336) $ \\
        $5$    & & & & $(2240, 280)$ & $ (5040, 560) $ & $ (10080, 1008) $ \\
  $ 6 $ &&&&& $    (12600, 1260) $ & $ (27720, 2520) $ \\
   \hline
\end{tabular}
\caption{Bidegree of the logarithmic determinant of a $d_1 \times d_2$ matrix}
\label{tab:logdet}
\end{table}

\begin{corollary} \label{cor:logdet}
Let $f_z(t)$ be a general bilinear form in $d_1+d_2$ variables, so
$X_A = \PP^{d_1-1} \times \PP^{d_2-1}$ is the
Segre variety of dimension $d = d_1+d_2-2$  and degree  $v = \binom{d}{d_1-1}$ in $\PP^{d_1d_2-1}$.
Then $\Lambda_A$ has bidegree $(dh,h)$, where $\, h = 2v + 2g - 2\,$ and
 $ \,g = 1 + (d_1 d_2-2d_1-2d_2+3)\,v /d$.
\end{corollary}

\begin{proof}
The Segre variety $X_A$ is smooth in $\PP^{d_1d_2-1}$.
Hence equality holds in Theorem~\ref{thm:fromHurwitz}.
The formula for the sectional genus $g$ of $X_A$ is found using the Adjunction Formula.
\end{proof}

We note that the Hurwitz form can also be used when the matroid of $A$ is not~connected.

\begin{example}
The toric threefold $X_A$ in Example \ref{ex:jointwoconics} is the join of two  plane conics.
It lives in $\PP^5$ and it is singular along those conics.
In particular, the singular locus of $X_A$ has codimension $2$. The sectional genus equals $g=1$,
so the Hurwitz degree is $h=8$. The Hurwitz form ${\rm Hu}_A$ is an irreducible polynomial
of degree $dh = 24$ in the $18$ entries of a $3 \times 6$ matrix $Y$.
The specialization in Theorem \ref{thm:fromHurwitz}
splits up into irreducible factors as follows:
\[ 
z_1^4 z_3^4 \, (4 z_1 z_3-z_2^2) \,  \Lambda_{A^{(1)}}  \,  z_4^4 z_6^4\, (4 z_4 z_6{-}z_5^2) \,
\,  \Lambda_{A^{(2)}}  \,  \nu_3^2 \,(s+\nu_3)^2.
\]
Here, $ \Lambda_{A^{(1)}}$ and $ \Lambda_{A^{(2)}}$ are 
the transformed logarithmic discriminants shown in 
(\ref{eq:twohypersurfaces}).
The factors depending only on $z$ correspond to two singular curves and four singular points. 
The factors $\nu_3^2(s+\nu_3)^2$ replace the factor $s^{(d-1)h}$ in Theorem \ref{thm:fromHurwitz} via the coordinate change~$Q$. 
\end{example}

\section{Leading Coefficients}
\label{sec6}
    
This section investigates initial forms of the logarithmic $A$-discriminant. In particular, we 
study the evaluation of $\Lambda_A$ when some $s$- or $\nu$-variables are set to zero. This further elucidates the relationship between the
$A$-discriminant $\Delta_A$ and the logarithmic $A$-discriminant $\Lambda_A$. We begin with the following fundamental fact about the
logarithmic $A$-discriminant.

\begin{theorem} \label{thm:DeltaAfactor}
The $A$-discriminant $\Delta_A$ is a factor in the evaluation of
$\Lambda_A$ at $s=0$.
\end{theorem}

\begin{proof}
We assume that $X_A$ is not dual defective. Otherwise $\Delta_A=1$ and there is nothing to prove.
By the Horn uniformization (Proposition \ref{prop:Horn}), 
the hypersurface $\{\Delta_A=0\}$ has a dense subset
whose points  are Hadamard products $x \star v$, where $x \in X_A$ and $v \in L^\perp=\ker(\hat A)$. We use the 
logarithmic Horn uniformization to parametrize a curve $\epsilon \mapsto \psi(\epsilon) $ inside $\Lambda_A$,
 whose limit $\lim_{\epsilon \to 0} \psi(\epsilon)$ 
has $z$-coordinates $x \star v$ and $s$-coordinate $0$.
It is our goal to show that every point $(x \star v, (0,-\nu))$
in $ \mathbb{P}^{n-1} \times \mathbb{P}^d$ arises as such a limit.
This will establish our assertion.

Our curve is defined as follows. For $w \in \mathbb{C}^d$, define $u_\epsilon = \mathbf{1} + \epsilon \, A^T w \in L$.
We set $\psi(\epsilon) = (x \star v \star u_\epsilon^{-1}, \hat{A} (v \star u_\epsilon^{-1})) \in \Lambda_A$. This uses the logarithmic Horn uniformization in Theorem \ref{thm:logpara}. Clearly, the limit for $\epsilon \to 0$ in the first coordinate equals $x \star v$ as desired. 

The point $\hat{A} (v \star u_\epsilon^{-1})$ in $\PP^d$ depends on the parameter vector $w$. Its first coordinate is 
\[ v \cdot u_\epsilon^{-1} \, = \, v \cdot ( {\bf 1} - \epsilon A^T w + O(\epsilon^2)) \, = \, O(\epsilon^2).\]
Here we use $v \in L^\perp$. The other entries are given by 
\[ A_v ( {\bf 1} - \epsilon A^T w + O(\epsilon^2)) \, = \,- \epsilon A_v A^T w + O(\epsilon^2).\] 
Hence, in $\PP^d$, we have $\lim_{\epsilon \to 0} \hat{A}(v \star u_\epsilon^{-1}) \, = \, (0, -A_vA^Tw)$.
In particular, $s = 0$, as desired.

 It remains to show that each point $(0, -\nu)\in \PP^d$ is such a limit. This rests on Lemma~\ref{lem:dualdefect},
 which says that
 the matrix $A_vA^T$ is invertible for generic $v$ if and only if $X_A$ is not dual defective.
We conclude that $A_vA^T$ has rank $d$, and hence $ (0,-\nu) = -\hat{A}_v \hat{A}^T w$
for some $w$.
\end{proof}

\begin{remark} Theorem  \ref{thm:DeltaAfactor} implies that $\Lambda_A$ is a hypersurface if $X_A$ is not dual defective.
Our proof highlights the connection between the
Horn uniformization and its logarithmic~cousin.
\end{remark}

Recall that $h = h(A)$ is the Hurwitz degree
of  $X_A$, 
and $\Lambda_A$ is of degree $h$ in $(s, \nu)$. Hence
\begin{equation} \label{eq:expandLambda} \Lambda_A \,\, = \,\, c_0(z) \cdot s^h\,+\, c_1(z) \cdot \nu_1^h \,+ \,
c_2(z) \cdot \nu_2^h \,+ \,\,\cdots \,\,+\,c_d(z) \cdot \nu_d^h \,+ \, 
\hbox{other terms}. \end{equation}
Our goal is to understand the $d+1$ leading coefficients $c_i \in \ZZ[z]$.
Theorem \ref{thm:DeltaAfactor} implies:
\begin{corollary}
The $A$-discriminant $\Delta_A$ divides the  leading coefficient $c_i(z)$ for all  $\,i \geq 1$.
\end{corollary}

We illustrate the goal with two low-dimensional examples.

\begin{example}[$d=1$] \label{eq:rnc}
Let $A = [0,1,\ldots,n\!-\!1]$, so $X_A$ is the rational normal curve, and
 $f_z(t) = \sum_{i=1}^{n} z_i t^{i-1}$ is a dense polynomial.
 Then $\Lambda_A(z,s,\nu)$
has bidegree $(2n\!-\!4,2n\!-\!4)$,~and
\begin{equation}
\label{eq:cfor1}
 c_1(z) \,=\, {\rm discr}_t (f_z(t)) \, = \, \Delta_A
\quad {\rm and} \quad c_0(z) \,=\, z_2^2 \cdot {\rm discr}_t (f'_z(t)). 
\end{equation}
If the toric curve $X_A$ is singular then
the degree of the discriminant $\Delta_A$ drops
and the same for $\Lambda_A$.
For instance, if $d=1$ and $n=3$ with
$A = [0,p,q]$ and $0 < p < q$ relatively prime, then
the logarithmic $A$-discriminant has bidegree $(q,q)$, and its leading coefficients are
$\,  c_1(z) \, = \, \Delta_A(z) \,=
\,(-q)^q  z_1^{q-p} z_3^p \,-\,  p^p (q-p)^{q-p} z_2^q \,\,$ and $\,\,
c_0(z)\,=\, z_2^q \,\,\,\hbox{times a constant}$.
\end{example}

\begin{example}[$d=2$] 
Let $f_z(t)$ be a general quadratic polynomial as in Example \ref{ex:tact1}. Its $A$-discriminant is
the determinant
$\Delta_A = 4 z_1 z_4 z_6-z_1 z_5^2-z_2^2 z_6+z_2 z_3 z_5-z_3^2 z_4$.
We compute
\begin{equation}
\label{c0c1c2A} \begin{matrix}
c_0(z) & = & z_2^2 \,z_3^2 \,(2 z_2 z_6-z_3 z_5)^2 \,(2 z_3 z_4 - z_2 z_5)^2,
\smallskip \\
c_1(z) & = & 4\, z_4\, \Delta_A\, (4 z_1 z_6-z_3^2) \,(z_1 z_5^2+z_2^2 z_6-z_2 z_3 z_5)^2,
\smallskip \\
c_2(z) & = & 4 \,z_6 \,\Delta_A\, (4 z_1 z_4-z_2^2) \,(z_1 z_5^2+z_3^2 z_4-z_2 z_3 z_5)^2.

\end{matrix}
\end{equation}
We here see the squares of four resultants, namely for the triples
$(x,\partial_1 f_z, \partial_2 f_z)$,
$(y,\partial_1 f_z, \partial_2 f_z)$,
$(x,\partial_1 f_z, f_z)$ and
$(y,\partial_2 f_z ,  f_z)$. 
This  will be explained by means of the following general result.
\end{example}

Let $\hat{A}_i \in \mathbb{Z}^{d \times n}$ be the matrix obtained by removing the $i$-th row from $\hat{A}$. 
Here we index the rows of $\hat{A}$ starting from zero, so $i$ runs from $0$ to $d$. In particular, we have $\hat{A}_0 = A$. 

\begin{proposition} \label{prop:leading}
Suppose that $X_A$ is smooth. For $i = 0, 1, \ldots, d$, up to a nonzero constant factor, the coefficient $c_i(z)$ is the evaluation of the Hurwitz form ${\rm Hu}_A$ at the matrix $Y = (\hat{A}_i)_z$. 
\end{proposition}
\begin{proof}
We know from Theorem \ref{thm:fromHurwitz} that the evaluation of the Hurwitz form ${\rm Hu}_A$ at the 
$d \times n$ matrix $T_{s, \nu} \hat{A}_z$ equals $C \cdot s^{h(d-1)} \Lambda_A$ for some $C \in \mathbb{Z} \backslash \{0\}$. From this, we see that 
\[ C \cdot c_0(z) \, = \, {\rm Hu}_A(T_{1, 0} \hat{A}_z) \, = \,  {\rm Hu}_A(A_z).\]
This proves the statement for $i = 0$. We now compute $c_i(z)$ for $i > 0$. Let $T_{s,\nu,i}$ be the matrix obtained by multiplying all but the $i$-th row of $T_{s, \nu}$ by $s^{-1}$. Since the Hurwitz form ${\rm Hu}_A(Y)$ is homogeneous of degree $h$ in each of the rows of $Y$, we have 
\[ {\rm Hu}_A(T_{s, \nu, i} \hat{A}_z) \, = \, s^{-h(d-1)} {\rm Hu}_A(T_{s, \nu} \hat{A}_z) \, = \, C \cdot \Lambda_A(z, s, \nu). \] 
Setting $s(\epsilon) = \epsilon  s$ and $\nu(\epsilon) = (\epsilon \nu_1, \ldots, \epsilon \nu_{i-1}, 1, \epsilon \nu_{i+1}, \ldots, \epsilon \nu_d)$, we have 
\[ C \cdot c_i(z) \, = \, \lim_{\epsilon \to 0} \, C \cdot \Lambda_A(z,s(\epsilon), \nu(\epsilon)) \, = \, \lim_{\epsilon \to 0} \, {\rm Hu}_A(T_{s(\epsilon), \nu(\epsilon), i} \hat{A}_z).\]
There exists a $d \times d$ matrix $Q$ with $\det Q = 1$ such that $\lim_{\epsilon \to 0} T_{s(\epsilon), \nu(\epsilon), i} \hat{A}_z = Q \cdot (\hat{A}_i)_z$. This proves the statement for $i > 0$, because ${\rm Hu}_A(Q (\hat{A}_i)_z) =  {\rm Hu}_A((\hat{A}_i)_z)$. 
\end{proof}

It remains to understand the specialized Hurwitz forms in Proposition \ref{prop:leading}.
This task~is interesting and challenging.
To tackle it, we write ${\rm md}_t$ for the operator 
that extracts the mixed discriminant \cite{CCDRS} of $d$ polynomials in $d$ variables.
Likewise, we write ${\rm res}_t$ for the resultant of $d+1$ polynomials in $d$ variables.
We begin with the relevant notions for $d=2$.

Fix polynomials  $l_1,l_2$ of degree $1$ and
$p_1,p_2$ of degree $k-1$.
These are dense polynomials  in $(t_1,t_2)$ with unknown coefficients.
Record degrees in these coefficients as vectors in $\NN^4$.
We consider the following specialization of the Hurwitz form for
the degree $k$ Veronese surface:
\begin{equation}
\label{eq:7factors} \begin{matrix} {\rm Hu}_A \bigl(\, l_1 \cdot p_1, \,l_2 \cdot p_2 \,\bigr)
\,\,\, = \,\,\,
{\rm md}_t( p_1,p_2) \cdot {\rm md}_t( l_1,p_2) \cdot 
{\rm md}_t( p_1,l_2) \cdot \qquad \qquad \qquad
\\ \qquad \cdot \,
{\rm res}_t(l_1,l_2,p_1)^2 \cdot
{\rm res}_t(l_1,l_2,p_2)^2 \cdot
{\rm res}_t(l_1,p_1,p_2)^2 \cdot
{\rm res}_t(l_2,p_1,p_2)^2.
\end{matrix} 
\end{equation}
This  is derived from 
\cite[Corollary 5.5]{BB}; see also
\cite[equation (7)]{Stu}. 
It is gratifying to see that the $\NN^4$-degrees in (\ref{eq:7factors}) add up correctly.
Indeed, using \cite[equation (1.4)]{CCDRS}, we check that
$$ \begin{small}  \begin{matrix} 
[3 k (k {-} 1), 3 k (k {-} 1), 3 k (k {-} 1), 3 k (k {-} 1)] \,\, =
\qquad \qquad \qquad \qquad \qquad \qquad \qquad \qquad 
\qquad \qquad \qquad \qquad \quad
\smallskip  \\
[0, 0, 3 (k{-}1) (k{-}2), 3 (k{-}1) (k{-}2)] \,+\,
[(k{-}1) (k{-}2), 0, 0, 2 k{-}4] \,+\,
[0, (k{-}1) (k{-}2), 2 k{-}4, 0]  \\
\,+ \, 2 \cdot [k{-}1, k{-}1, 1, 0]
\,+ \,2 \cdot [k{-}1, k{-}1, 0, 1]
\, + \,2\cdot [(k{-}1)^2, 0, k{-}1, k{-}1]
\, + \,2 \cdot [0, (k{-}1)^2, k{-}1, k{-}1].  \end{matrix}
\end{small}
$$
Furthermore, if $q$ is a general irreducible polynomial of degree $k$ in $(t_1,t_2)$ then we also have
\begin{equation}
\label{eq:3factors}  {\rm Hu}_A \bigl(\, q \,, \,l_2 \cdot p_2 \,\bigr) \,\,\,=\,\,\,
{\rm md}_t(q,l_2) \cdot {\rm md}_t(q,p_2) \cdot
{\rm res}_t(q ,l_2,p_2)^2.
\end{equation}

Our formulas describe degenerations of
Hurwitz forms. This topic appeared in \cite[Section 4]{Stu}
and it was further developed in \cite[Section 5]{BB}.
The factors are Hurwitz forms of 
top strata and squares of Chow forms of codimension one strata.
This explains the squares in (\ref{c0c1c2A}).
Indeed, (\ref{eq:7factors}) and (\ref{eq:3factors}) apply directly to the
scenarios that arise in Proposition \ref{prop:leading}.

\begin{corollary}[$d=2$]
Let $A \in \mathbb{Z}^{2 \times \binom{k+2}{2}}$ be the support of a generic bivariate polynomial $f_z$ of degree $k$. The coefficient $c_0(z)$ in (\ref{eq:expandLambda}) equals (\ref{eq:7factors}) 
with $l_1 = t_1$, $l_2 = t_2$,
$p_1 = \partial_1 f_z$ and $p_2 = \partial_2 f_z$.
The other coefficients $c_i(z)$ are
(\ref{eq:3factors}) with $q = f_z$,
$l_2 = t_j$ and $p_2 = \partial_j f_z$ for $j \not= i$.
\end{corollary}

We note that ${\rm md}_t ( f_z, t_2) = 4z_1z_4-z_2^2$ and
${\rm md}_t(f_z,t_1) = 4 z_1 z_6 - z_3^2$ in (\ref{c0c1c2A}).
Theorem~\ref{thm:DeltaAfactor} shows itself in the fact that $\Delta_A$ divides both
${\rm md}_t( f_z, \partial_1 f_z) = z_6 \Delta_A$ and  ${\rm md}_t( f_z, \partial_2 f_z) =
z_4 \Delta_A$.
The divisibility is usually strict.
We now state the same fact for all $d$ and all supports $A$.

\begin{proposition} \label{prop:divides}
Let $f_z(t)$ be the  general polynomial in $t=(t_1,\ldots,t_d)$ with support $A$.
Fix $i \in \{1,2,\ldots,d\}$ and consider the following specialization of the mixed discriminant:
\begin{equation}
\label{eq:mdt}
 {\rm md}_t \bigl( \,f_z,\, \partial_1 f_z,\ldots,\partial_{i-1} f_z, \,\partial_{i+1} f_z, \ldots,\partial_d f_z\, \bigr).
 \end{equation}
Here ${\rm md}_t$ is the mixed discriminant with respect to the monomial supports of $f_z$ and the $d-1$ derivatives appearing in \eqref{eq:mdt}. If ${\rm md}_t \neq 1$, then it has
 the $A$-discriminant $\Delta_A$ as a~factor.
 \end{proposition}

\begin{proof}
Without loss of generality, we prove the proposition for $i = d$. Suppose that $p \in (\mathbb{C}^*)^d$ is a singular point of the hypersurface defined by $f_z$, $z \in \nabla_A$. This means $f_z(p) = 0$ and $\partial_i f_z (p) = 0$ for $i = 1, \ldots d$. The following deformation of $f_z$ preserves the monomial support: 
\[ f_{z, \epsilon}(t) \, = \, f_z(t) + \epsilon \, \sum_{i = 1}^{d-1} t_i \, \partial_i f_z(t). \]
We will show that ${\rm md}_t(f_{z, \epsilon}, \partial_1 f_z, \ldots, \partial_{d-1} f_z) = 0$ for any $\epsilon \neq 0$, which implies by continuity that ${\rm md}_t(f_{z}, \partial_1 f_z, \ldots, \partial_{d-1} f_z) = 0$. This proves that each $z \in \nabla_A$ is contained in the mixed discriminant hypersurface, so $\Delta_A$ must factor the mixed discriminant polynomial. 

Consider the Cayley polynomial $g_\epsilon(t, y) = f_{z, \epsilon}(t) + y_1 \, \partial_1 f_z(t) + \cdots + y_{d-1} \, \partial_{d-1} f_{z}(t)$. Let $\hat{p}$ be the point $(p, -\epsilon \,  p_1, \ldots, -\epsilon \,  p_{d-1}) \in (\mathbb{C}^*)^{2d-1}$. One checks that $g_\epsilon$ is singular at $\hat{p}$: 
\[ g_\epsilon(\hat{p}) \, = \,  \partial_{t_1} g_\epsilon (\hat{p}) \, = \, \cdots \, = \,  \partial_{t_d} g_\epsilon (\hat{p}) \, = \, \partial_{y_1} g_\epsilon (\hat{p}) \, = \, \cdots \, = \,  \partial_{y_{d-1}} g_\epsilon (\hat{p}) \, = \, 0. \]
Therefore, $g_\epsilon$ belongs to the discriminant variety $\{\Delta_{\rm Cay} = 0\}$ associated with the Cayley configuration of the monomial supports in \eqref{eq:mdt}. That Cayley discriminant variety is indeed a hypersurface, and it equals $\{{\rm md}_t = 0\}$ by \cite[Theorem 2.1]{CCDRS}. This uses our assumption that ${\rm md}_t \neq 1$. We conclude that $\Delta_{\rm Cay}(g_\epsilon) \, = \, {\rm md}_t(f_{z, \epsilon}, \partial_1 f_z, \ldots, \partial_{d-1} f_z) = 0$, as desired. 
\end{proof}

% \begin{proposition} \label{prop:divides}
% Let $f_z(t)$ be the 
% general polynomial in $t=(t_1,\ldots,t_d)$
% with support $A$. Then, for every $i \in \{1,2,\ldots,d\}$,
% the $A$-discriminant $\Delta_A$ divides the mixed discriminant
% \begin{equation}
% \label{eq:mdt}
%  {\rm md}_t \bigl( f_z, \partial_1 f_z,\ldots,\partial_{i-1} f_z, \partial_{i+1} f_z, \ldots,\partial_d f_z \bigr). 
%  \end{equation}
% \end{proposition}

% \begin{proof}
% At a generic point $z^* \in \nabla_A$, the hypersurface $\{ f = 0\} \subset (\mathbb{C}^*)^d$ has a singular point. That is, there exists $t^* \in (\mathbb{C}^*)^d$ such that $f(t^*) = \partial_1 f(t^*) = \cdots = \partial_d f(t^*) = 0$. Such a singular point is also a common zero of the $d$ polynomials in $d$ variables in \eqref{eq:mdt}.
% Moreover, the Jacobian matrix evaluated at the common zero $t^*$ is singular, because the row corresponding to $f_z$ is zero.
% Hence, the Jacobian determinant vanishes at any such singular point, and the specialized mixed discriminant must vanish at the coefficients $z^*$. 
% \end{proof}

The factorizations in (\ref{eq:7factors})
and (\ref{eq:3factors}) can be generalized to $d \geq 3$, and so can Proposition \ref{prop:divides}.
The resulting combinatorial structures are rich and they deserve to be studied further.

\smallskip

We now shift gears and briefly discuss a different class of leading terms.
These were studied recently in the physics literature, namely in \cite{CCG}.
Recall from Corollary \ref{cor:ER} that $E_A(z)\cdot R_A(s,-\nu)$
vanishes when one of the $v$ critical points
escapes from the very affine variety $(\CC^*)^d\setminus V(f_z)$. 
It is natural to ask a more refined
question, namely which factor of $E_A(z)R_A(s,-\nu)$
vanishes when a particular coordinate
$t_i$ tends to zero or to infinity.

Let $K = \QQ(z,s,\nu)$ be the rational function field
in our parameters, and consider the ideal
$$ I \,\,:= \,\,\biggl\langle  \nu_1 f_z + s t_1 \frac{\partial f_z}{\partial t_1} ,\,\ldots\,,\,
\nu_d f_z + s t_d \frac{\partial f}{\partial t_d}, \,t_0 t_1 \cdots t_d f_z -1 \biggr\rangle
\,\,\,\subset \,\,\, K[t_0,t_1,t_2,\ldots,t_d].
$$
Here $t_0$ is a new variable. The last equation specifies $\mathcal{U}_z$ as a variety
in affine $(d+1)$-space. The ideal $I$ is zero-dimensional and radical. 
Its variety is the critical locus, so its degree is~$v$.

The elimination ideal that describes the $i$-th coordinate for each critical point is principal:
\begin{equation}
\label{eq:principal}
  \quad  I \,\,\cap \,\, K[t_i] \,\, = \,\, \bigl\langle p_i(t_i) \bigr\rangle  \qquad
{\rm for} \,\,\, i =0,1,\ldots,d. 
\end{equation}
The univariate polynomial $p_i(t_i)$ is irreducible and we  expect its degree to be
 $v={\rm vol}(P)$. It is unique up to units, so we assume its coefficients to be
 relatively prime in the polynomial ring $\QQ[z,s,\nu]$.
 Under that hypothesis, let $c_i^+(z,s,\nu)$ denote the leading coefficient of $p_i(t_i)$,
and let $c_i^{-}(z,s,\nu)$ denote the trailing coefficient.
These are polynomials in $n+1+d$ variables with coefficients in the field
$\QQ$ of rational numbers. Their factors appear in $\,E_A(z) \cdot R_A(s,-\nu)$.

\begin{example}[Triangular Prism]  \label{ex:toblerone}
Fix $A$ as in Example \ref{ex:23independence}. The coordinates of the 
$v=3$ critical points are the zeros of the univariate polynomials $p_1(t_1),p_2(t_2),p_3(t_3)$.
Here we find
$$ 
\begin{matrix} c_1^+ \,=\,  z_2 z_5 (z_2 z_6-z_3 z_5)\cdot (s+\nu_1+\nu_2)^2,
             &       & c_1^- \,=\,  z_1 z_4(z_1 z_6-z_3 z_4) \cdot \nu_1^2,  \\
 c_2^+ \,=\,     z_3 z_6 (z_2 z_6-z_3 z_5)\cdot (s+\nu_1+\nu_2)^2,& & 
 c_2^- \,=\,   (z_2z_4-z_1z_5)z_1 z_4 \cdot \nu_2^2 , \\            
c_3^+ \,=\,     z_4 z_5 z_6 \cdot (s+\nu_3) , \qquad \qquad \qquad \quad \,
& &  c_3^- \,=\, z_1 z_2 z_3 \cdot \nu_3.\qquad \qquad \quad \, \\
\end{matrix}
$$
The non-vanishing of these expressions ensures that the three critical points remain in $(\CC^*)^3$.
\end{example}

We leave it to a future project to determine the  factors of $c_i^+(z,s,\nu)$ and $c_i^{-}(z,s,\nu)$
for arbitrary $A$. These are the
{\em Leviathans} from \cite{CCG} in the GKZ setting.
The answer will rest on
 which faces of
 Newton polytope $P = {\rm conv}(A)$ are visible in the 
 coordinate directions $\pm e_i$.

\section{Arrangements of Hypersurfaces}
\label{sec7}

In this section we consider an arrangement of
$m \geq 2$ hypersurfaces in the torus $(\CC^*)^d$,
specified by a tuple $\boldA$ of configurations
$A_1,\ldots,A_m$ in $\ZZ^d$. The $i$th hypersurface is defined~by
$$ \quad  f_i(t) \,\, = \,\, \sum_{j=1}^{n_i} z_{ij} \cdot t^{a_{j}^{(i)}} \qquad {\rm for} \,\,\, i=1,2,\ldots,m. $$
Here, $A_i = (a_{1}^{(i)},a_{2}^{(i)},\ldots,a_{n_i}^{(i)}) \in \ZZ^{d \times n_i}$. 
Our $m$ polynomials have  $n = n_1\!+\!n_2\!+\!\cdots \!+\! n_m$ coefficients in total.
These are $n$ distinct variables,  which we collect in a vector of length $n$:
\begin{equation} \label{eq:zm} z \,=\, (z_1,z_2,\ldots,z_m) \,=\,
(z_{11}, \ldots, z_{1n_1}, z_{21}, \ldots, z_{2n_2}, \ldots, z_{m1}, \ldots, z_{mn_m}). \end{equation}
The subvector $z_i = (z_{i1},\ldots,z_{in_i})$
gives coordinates on $\PP^{n_i-1}$, so the vector $z$ 
gives coordinates on the product of projective spaces
$\PP^{n_1-1} \times \PP^{n_2-1} \times \cdots \times \PP^{n_m-1}$.
Suppressing the subscript $z$ for convenience,
as in (\ref{eq:arrangement}),
our Euler integral (\ref{eq:eulerint}) now takes the more general form

\begin{equation} \label{eq:eulerintm}
 \int_\Gamma f(t)^s t^\nu {\rm d} t \,\,=\,\,
  \int_\Gamma f_1(t)^{s_1} \cdots f_m(t)^{s_m}
\, t_1^{\nu_1} \cdots t_d^{\nu_d} \,{\rm d} t_1 \cdots {\rm d} t_d.
 \end{equation}

As before in (\ref{eq:loglike}), the log-likelihood function is the logarithm of the integrand:
\begin{equation}
\label{eq:loglikem}
 {\rm log}\bigl( f(t)^s t^\nu \bigr) \,\, = \,\, 
\sum_{i=1}^m s_i \cdot {\rm log}( f_i(t))\,+\,
\sum_{j=1}^d \nu_j \cdot {\rm log}(t_j). 
\end{equation}
We take $(s,\nu) = (s_1,\ldots,s_m,\nu_1,\ldots,\nu_d)$
as homogeneous coordinates for the space $\PP^{m+d-1}$.
All our results in Sections \ref{sec2}-\ref{sec6} were for $m = 1$,
and we now extend them to $m \geq 2$. In fact, we 
shall establish a strong resemblance to the $m = 1$ case, which allows us to omit~proofs. 

We will see below that, after clearing denominators and introducing new auxiliary variables $y_i = s_if_i^{-1}, i = 1, \ldots, m$, all partial derivatives of (\ref{eq:loglikem})
have the same support, namely
$$ \Cay(\boldA) \,\,= \,\,\begin{small} \begin{bmatrix}
{\bf 1} & {\bf 0} & \cdots & {\bf 0} \\
{\bf 0} & {\bf 1} & \cdots & {\bf 0} \\
 \vdots & \vdots & \ddots & \vdots \\
{\bf 0} & {\bf 0} & \cdots & {\bf 1} \\
A_1 & A_2 & \cdots & A_m 
\end{bmatrix}. \end{small}
$$
This is known as the {\em Cayley configuration} of the tuple $\boldA = (A_1,\ldots,A_m)$. 
In the special case when
 $d=m$ and $A = A_1 = \cdots = A_d$, we recover the matrix
 $\Cay(\boldA) = A^{(d)}$ seen in (\ref{eq:A(d)}). We assume that the matroid of $\Cay(\boldA)$ is connected, generalizing our standing assumption~on~$\hat{A}$.

We additionally assume that the columns of ${\rm Cay}(\boldA)$, after dropping the first row, affinely generate the lattice $\mathbb{Z}^{d+m-1}$. Then the normalized volume of the Cayley polytope ${\rm conv}(\Cay(\boldA))$ is
the number $v$ of critical points of the function (\ref{eq:loglikem}).
This is made precise in the second part of Theorem \ref{thm:ERLm}.
The critical points are algebraic functions of $z,s,\nu$. The {\em logarithmic
$\boldA$-discriminant $\Lambda_{\boldA}$}
 vanishes when two
of these critical points come together.

The hypersurface defined by $\Lambda_{\boldA}$ lives naturally in the product of projective spaces
$$ \PP^{n_1-1} \times \PP^{n_2-1} \times \cdots \times \PP^{n_m-1} \,\times\, \PP^{m+d-1}. $$
The degree of $\Lambda_{\boldA}$ is a vector in $\NN^{m+1}$. Its first
coordinate is the degree of $\Lambda_{\boldA}$ in the coefficient vector $z_1$ of the first
polynomial $f_1(t)$. Its last coordinate is the degree of $\Lambda_{\boldA}$ in $(s,\nu)$.

\begin{example}[Two lines] \label{ex:twolines}
Let $d=m=2$ and $A_1 = A_2 = \begin{footnotesize}
\begin{bmatrix} 0 & \! 1\! & 0 \\ 0 & \!0\! & 1 \end{bmatrix} \end{footnotesize}$. 
In (\ref{eq:loglikem}) we consider
$$ (t_1,t_2) \,\mapsto \,  s_1 \cdot {\rm log}( z_{11}\!+\! z_{12} t_1 \!+\! z_{13} t_2 )
+ s_2 \cdot {\rm log}(z_{21} + z_{22} t_1 + z_{23} t_2)
+ \nu_1 {\rm log}(t_1) + \nu_2 {\rm log}(t_2).
$$
This has $v=3$ critical points because 
${\rm conv}({\rm Cay}(\boldA))$ is the triangular prism in Example~\ref{ex:23independence}.
The critical points are algebraic functions in
the $3+3+4$ unknowns $z_1,z_2,(s,\nu)$.
The logarithmic $\boldA$-discriminant $\Lambda_{\boldA}$ is a sum of
$305$ monomials of degree $(6,6,4)$. It defines
 an irreducible hypersurface  in $\PP^2 \times \PP^2 \times  \PP^3$.
The degree in $(s,\nu)$ is $4$. This matches
\cite[Theorem 1.2]{KKT}. The complexity increases for
$m=3$, when our arrangement is a triple
of lines in $(\CC^*)^2$. The logarithmic ${\boldA}$-discriminant $\Lambda_{\boldA}$
is a hypersurface of degree $(16,16,16,12)$ in $(\PP^2)^3 \times \PP^4$.
\end{example}

As in our earlier discussion for $m=1$, we seek
 critical points in  the very affine~variety 
$$\mathcal{U}_z \,\,=\,\, (\mathbb{C}^*)^d \setminus V(f_1 f_2 \cdots f_m). $$
The Euler characteristic of $\mathcal{U}_z$
 is $(-1)^d v$, where $v = {\rm vol}({\rm Cay}(\boldA))$, 
provided $z$ is not~a zero of
the principal $\Cay(\boldA)$-determinant
$E_{\boldA}(z):=E_{\Cay(\boldA)}(z)$.
This follows from \cite[Theorem~2.2]{FM}.

\begin{example}[Conics in the plane]
It is instructive to compare Examples \ref{ex:tact1} and \ref{ex:twolines}.
If $m=1$ and $f(t)$ is quadratic
in $(t_1,t_2)$ then $\mathcal{U}_z = ( \CC^*)^2 \backslash V(f)$ has
Euler characteristic $v=4$, so $f(t)^s t^\nu$ has $4$ critical points.
But if $f = f_1 f_2$ factors into two linear factors,
then the Euler characteristic drops to $3<v$. Now,  $f(t)^s t^\nu$ has
only $3$ critical points. Let us specialize $f=f_1f_2$ in
the discriminant $\Lambda_A$ from Example \ref{ex:tact1}.
The result has degree $(12,12,6)$. It has an extraneous factor of
degree $(6,6,2)$, namely
$(z_{11} z_{22} - z_{12} z_{21})^2 \cdot
  (z_{11} z_{23} - z_{13} z_{21})^2 \cdot
  (z_{12} z_{23} - z_{13} z_{22})^2 \cdot s^2$.
  The other factor has $233$ terms of degree $(6,6,4)$.
  It equals $\Lambda_\boldA(z,s,s,\nu_1,\nu_2)$.
\end{example}

The critical points are defined by the partial derivatives of \eqref{eq:loglikem}:
\begin{equation}\label{eq:criteqsm} \quad
 \frac{\nu_i}{t_i} \,+\, \sum_{j = 1}^m s_j \, \frac{ \partial_i \,f_j(t)}{ f_j(t)}
    \, \,= \,\, 0 \qquad {\rm for} \quad  i = 1, \ldots, d.
\end{equation}
The Cayley configuration enters the story by introducing new variables $y=(y_1,\dots,y_m)$:
$$ y_k f_k(t) \, = \,  s_k \,\,\,\, {\rm for} \,\,\, k=1,\ldots,m \quad {\rm and}
 \quad \nu_i + \sum_{j=1}^m y_j \, t_i \, 
 \partial_i f_j(t)     \, = \, 0 \,\,\,\,{\rm for} \,\, \, i = 1, \ldots, d.
$$
Setting $x = (x^{(1)}, x^{(2)}, \ldots, x^{(m)})$ with $x^{(j)} = (y_j t^{a^{(j)}_1}, \ldots, y_j t^{a^{(j)}_{n_j}})$, the critical equations are 
\begin{equation}\label{eq:criteqm}
  x\in X_{\Cay(\boldA)} \quad\text{s.t.}\quad {\rm rank} \begin{bmatrix} \Cay(\boldA)\diag(z)\,x
   &  \mu \end{bmatrix} \, \leq \, 1,
\end{equation}
where $\mu = (s, -\nu)$. This defines a subscheme ${\rm Crit}$ of the projective toric variety $X_{{\rm Cay}(\boldA)}$ as in \eqref{eq:Crit}. Let $Q \in \mathbb{Z}^{(d+m) \times (d+m)}$ be the row operation which replaces the first row of ${\rm Cay}(\boldA)$ by the sum of its first $m$ rows. The first row of $Q {\rm Cay}(\boldA)$ is the all-ones vector, and $Q \mu = (s_1 + \cdots + s_m,s_2, \ldots, s_m, -\nu_1, \ldots, -\nu_d)$. Applying $Q$ to the left of $\begin{bmatrix} \Cay(\boldA)\diag(z)\,x
   &  \mu \end{bmatrix}$, \eqref{eq:criteqm} become the critical equations for the following integrand with $d + m - 1$ variables:
   \begin{equation}  \label{eq:cayleyintegrand} g(t,y)^{s_1 + \cdots + s_m} \, y_2^{-s_2} \cdots y_m^{-s_m} \, t_1^{\nu_1} \cdots t_d^{\nu_d}, \quad \text{where} \quad g(t,y) \, = \, f_1(t) + \sum_{j = 2}^m y_j \, f_j(t). \end{equation}
A critical point in $(\mathbb{C}^*)^{d + m-1} \backslash V(g)$ corresponds to a solution of \eqref{eq:criteqsm} in ${\cal U}_z$ by projecting to the $t$-coordinates. Our next result characterizes the inverse  $t \mapsto (t,y)$ of this correspondence.

\begin{lemma} \label{lem:japan}
If $s_1 s_2 \cdots s_m (s_1 + \cdots + s_m) \neq 0$, then the critical points of \eqref{eq:loglikem} are in
bijection with those of the logarithm of \eqref{eq:cayleyintegrand}. 
The latter are recovered from the former by the formula
\begin{equation}
\label{eq:corrinverse}
\quad  y_i \, = \, \frac{s_i \, f_1(t)}{s_1 \, f_i(t)} \qquad {\rm for} \,\,\, i = 2, \ldots, m. 
 \end{equation}
\end{lemma}

\begin{proof}
Let $s_+ = s_1 + \cdots +s_m$.
The condition $s_1 s_2\cdots s_m \,s_+ \neq 0$ is found by computing the inverse of the map $(t, y) \mapsto t$ from the critical points of  \eqref{eq:cayleyintegrand} to the solutions of \eqref{eq:criteqsm}.
Setting the $y$-derivatives of the logarithm of \eqref{eq:cayleyintegrand} to zero gives $y_j = s_jg(t,y)s_+^{-1}f_j(t)^{-1}$.
 This yields $m-1$ affine-linear conditions on $y_2, \ldots, y_m$. After homogenizing $(y_1 = 1)$, these are given by
\begin{equation} \label{eq:auxlineariny}\sum_{j=1}^{m} y_j \, f_j(t) \,-\, \frac{s_+ f_i(t)}{s_i} \, y_i \,\, =\, \, 0 
 \qquad {\rm for} \,\,\, i = 2, \ldots, m.  \end{equation}
We claim that  the $(m-1) \times m$ matrix for these homogeneous linear equations has full rank $m-1$.
 The $(m-1)\times(m-1)$ minor indexed by $y_2, \ldots, y_m$ is the determinant of a diagonal matrix plus a rank-one matrix. 
 By the Matrix Determinant Lemma, this minor equals
 \[
 \prod_{i=2}^m\left(-\frac{s_+}{s_i}f_i(t)\right)\left(1-\frac{\sum_{i=2}^m s_i}{s_+}\right)\, \,
 = \,\, \frac{s_1}{s_+}\prod_{i=2}^m\left(-\frac{s_+}{s_i}f_i(t)\right)\,\,\neq \,\,0.
\]
Hence, \eqref{eq:auxlineariny} has a unique solution with $y_1 = 1$. Also, using \eqref{eq:auxlineariny} we find
\[ \sum_{j = 1}^m y_j \, f_j \, = \, y_1 \, f_1 + \sum_{i = 2}^m y_i \, f_i \, = \, y_1 \, f_1 + \big( \sum_{j=1}^m y_j \, f_j)\sum_{i = 2}^m \frac{s_i}{s_+}, \quad \text{so} \quad \sum_{j} y_j \, f_j \, = \, s_1^{-1} s_+ f_1\, y_1.\]
The map (\ref{eq:corrinverse}) is found by substituting the righthand side of the latter equation into \eqref{eq:auxlineariny}.\end{proof}

Let $\widetilde{{\rm Cay}(\boldA)}$ be the $(d + m -1) \times n$ matrix obtained by deleting the first row of ${\rm Cay}(\boldA)$.

\begin{definition}
The {\em logarithmic $\boldA$-discriminant $\Lambda_{\boldA}$}
 is the evaluation of $\Lambda_{\widetilde{{\rm Cay}(\boldA)}}$ at the coefficients $z$ 
 in \eqref{eq:zm} and the exponents $Q \mu = (s_1 + \cdots + s_m,s_2, \ldots, s_m, -\nu_1, \ldots, -\nu_d)$.
\end{definition}

By construction, (the zero locus of) $\Lambda_{\boldA}$ is the closure of the set of parameters $(z,(s,\nu))$ for which  $s_1 \cdots s_m (s_1 + s_2 + \cdots + s_m) \neq 0$ and (\ref{eq:loglikem})
has a degenerate critical point in $\mathcal{U}_z$.
   
We now present the extension of Theorems~\ref{thm:irredcodim1} and
\ref{thm:ERL} to the arrangement case $m\geq 2$. The resonance arrangement $R_{\boldA}(s,-\nu)$ is the product 
of the linear forms defining the
 facets of the convex 
 polyhedral cone in $ \mathbb{R}^{d+m}$ that is generated by the columns of the matrix ${\rm Cay}(\boldA)$. 
Building on our earlier considerations,
the proof of Theorem \ref{thm:ERLm} follows from Lemma \ref{lem:japan}.

\begin{theorem}\label{thm:ERLm}
The logarithmic $\boldA$-discriminant $\Lambda_{\boldA}$ is an irreducible
hypersurface in 
$ \PP^{n_1-1} \times \PP^{n_2-1} \times \cdots \times \PP^{n_m-1} \,\times\, \PP^{m+d-1}$.
Suppose that the parameters satisfy $s_1 \cdots s_m (s_1 + \cdots + s_m) \neq 0$. Then the log-likelihood function (\ref{eq:loglikem}) has
  precisely $v$ distinct critical points in $\mathcal{U}_z$ if and only~if
 \begin{equation}\label{eq:tripleproductm}
    E_{\boldA}(z)\cdot R_{\boldA}(s,-\nu)\cdot\Lambda_{\boldA}(z,s,\nu)\,\neq \,0.
  \end{equation}
  Here   $R_{\boldA}:= R_{\Cay(\boldA)}$ is the resonance arrangement of $\Cay(\boldA)$ in $\PP^{d+m-1}$.
\end{theorem}

The equations \eqref{eq:criteqm} define
a linear space of dimension $n-m-d$ in $\PP^{n-1}$
which intersects $X_{\Cay(\boldA)}$ in $v$ points
when the parameters are generic. That linear space can be written as in \eqref{eq:coeffmtx}, where $T_{s,\nu}$ is now a $(d+m-1)\times(d+m)$ matrix whose first column is $(-s_2,\dots,-s_m,\nu_1,\dots,\nu_d)$, and whose remaining columns form the diagonal matrix $\diag(s_1,\dots,s_1)$. 
To identify special parameters, we substitute
our linear space into the Hurwitz form of 
$X_{\Cay(\boldA)}$. Namely,
the coefficients of our linear equations form
an $(m+d-1) \times n$ matrix $Y$. These
are the primal Stiefel coordinates. The first sentence in Theorem  \ref{thm:fromHurwitz}
holds verbatim for $m \geq 2$.
The degree statement generalizes too, but
we will not pursue this here.

\begin{example}[$d=m=2$] \label{ex:23independencec}
We revisit the two lines in Example~\ref{ex:twolines}.
The Cayley polytope equals the triangular prism 
in Example \ref{ex:23independence}, but now with a
different matrix representation:
\begin{equation}
\label{eq:newtoblerone} {\rm Cay}(\boldA) \,\, = \,\,
\begin{small} \begin{bmatrix}
 \,1 & 1 & 1 \, &\, 0 & 0 & 0 \,\,\\
 \,0 & 0 & 0 \, &\, 1 & 1 & 1\, \,\\
 \,0 & 1 & 0 \, &\, 0 & 1 & 0\, \,\\ 
 \,0 & 0 & 1 \, &\, 0 & 0 & 1 \,\,\\
 \end{bmatrix}. \end{small}
\end{equation}
We intersect the toric threefold $\,X_{\Cay(\boldA)} \simeq \PP^1 \times \PP^2\,$
 with three hyperplanes in $\,\PP^5$, given by
$$ 
Y \,\, = \,\,
\begin{bmatrix}
-z_{11}s_2  &  -z_{12}s_2  &  -z_{13}s_2  &  z_{21}s_1  &  z_{22}s_1  &  z_{23}s_1\\
z_{11}\nu_ 1  &  z_{12}(s_1+\nu_1) &  z_{13}\nu_1 &  0 &  z_{22}s_1 &  0 \\
z_{11}\nu_2 &  z_{12}\nu_2 &  z_{13}(s_1+\nu_2) &  0 &  0 &  z_{23}s_1 \\
\end{bmatrix}.
$$
We plug $Y$ into  the Hurwitz form of $X_{\Cay(\boldA)}$,
which is the $2 \times 3 \times 3$ hyperdeterminant, by Example \ref{ex:23independenceb}.
The result equals $s_1^8$ times
 $\Lambda_{\boldA}$, which has $305$ monomials, by Example~\ref{ex:twolines}.
\end{example}

The logarithmic Horn uniformization
(Theorem~\ref{thm:logpara}) also generalizes from $\hat A$ to $\Cay(\boldA)$.

\begin{theorem} \label{thm:hornm2}
  Let $L\subset\RR^n$ be the row space of $\Cay(\boldA)$ and $L^\perp\subset\RR^{n}$ the kernel of $\Cay(\boldA)$.   
  The hypersurface defined by the logarithmic $\boldA$-discriminant $\Lambda_{\boldA}$ is the image of the map 
  \begin{equation}
  \label{eq:asin2m}
  X_{{\rm Cay}(\boldA)} \times (L^\perp \star L^{-1})  \, \dashrightarrow\, \PP^{n_1-1}\times\cdots\times\PP^{n_m-1} \times \PP^{d+m-1}\,,\,\,\,
  (x,\phi)\, \mapsto \,\bigl(\,x \star \phi , \,\Cay(\boldA) \phi\,\bigr). 
  \end{equation}
\end{theorem}

The proof of this theorem is analogous to that in Section \ref{sec3}.
The formula in (\ref{eq:asin2m}) will be made explicit in the next section,
where it serves as the foundation of our numerical method.

\section{Computational Methods}
\label{sec8}

The logarithmic $\boldA$-discriminant $\Lambda_\boldA$ can be computed via elimination.
For simplicity, we  describe this now for the $m =1 $ version $\Lambda_A$
from Sections \ref{sec2}--\ref{sec6}.
The input is the $d \times n$ matrix~$A$,
which we then turn into a polynomial $f_z(t)$
with unknown coefficients $z = (z_1,\ldots,z_n)$.
The critical equations are given in \eqref{eq:criteqs}.
The existence of degenerate critical points is
  expressed by the additional equation $\det (\hat{A}_{z \star \Phi_A(t)} \hat{A}^T) = 0$ 
  which defines the ramification locus ${\cal R}$ of $\Psi$. This equation was found in the proof of Theorem \ref{thm:logpara}, and $\Phi_A$ is the monomial map in \eqref{eq:PhiA}.
  
  We perform symbolic computations in the polynomial ring
  $\QQ[t,y, z,s,\nu]$, which comprises $ d+2+n+(1+d) $ variables.
  Our incidence correspondence is now encoded by the ideal
  \begin{equation}\label{eq:idealtoelim}
  I\,\,=\,\,
  \bigl\langle \, \nu_1 + s \, y_1 \, t_1 \partial_1 \, f_z, \ldots, \nu_d + s \, y_1 \,  t_d \, \partial_d \, f_z , \, \det (\hat{A}_{z \star \Phi_A(t)} \hat{A}^T),  \, f_z \, y_1 -1, \, t_1 \cdots t_d \cdot y_2 -1\bigr\rangle.
\end{equation}
The last two generators ensure that the solutions lie in $\mathcal{U}_z$. Let $E = I \cap \mathbb{Q}[z, s, \nu]$ be the elimination ideal obtained by eliminating $t$ and $y$. 
Under the assumptions in the Introduction, the ideal $E$
 is principal, and the desired logarithmic $A$-discriminant $\Lambda_A$ is its generator. 

The elimination step can, in principle, be carried out with
any computer algebra system.
The following snippet shows how to use the system {\tt Oscar}
\cite{oscarbook} within {\tt Julia} for this computation. In this snippet, 
we assume that the entries of $A$ are nonnegative integers:
\begin{oscarcode}
using Oscar
function get_log_discr_elim(A)
    d,n = size(A)
    R,t,y,z,s,v = polynomial_ring(QQ,:t=>1:d,:y=>1:2,:z=>1:n,:s=>1:1,:v =>1:d)
    Ahat = matrix_space(ZZ,d+1,n)([ones(Int,1,n); A])
    p = [z[i]*prod(t.^A[:,i]) for i = 1:n]
    H = det(Ahat*diagonal_matrix(p)*transpose(Ahat))
    f = sum([z[i]*prod(t.^A[:,i]) for i = 1:n])
    I = ideal([v+(A*p).*(s[1]*y[1]); H; f*y[1]-1; prod(t)*y[2]-1])
    E = eliminate(I,[t;y]) 
end
\end{oscarcode}
Here the auxiliary variable \texttt{p} encodes the Hadamard product $z \star \Phi_A(t)$, and the toric derivatives of $f_z$ are conveniently written as the matrix-vector product $A (z \star \Phi_A(t))$. With the function \texttt{get\_log\_discr\_elim}, the logarithmic $A$-discriminant in Example~\ref{ex:square} is found via 

\begin{oscarcode}
  A = [0 1 0 1; 0 0 1 1]
  get_log_discr_elim(A)
\end{oscarcode}

Another approach is to use Theorem~\ref{thm:fromHurwitz} to compute $\Lambda_A$.
This can be preferable for configurations $A$ for which it is known how to
compute the  Hurwitz form of $A$ effectively. This scenario arises when $A$ is a product of simplices,
like the triangular prism $\Delta_1 \times \Delta_2$ in Examples \ref{ex:23independence},
 \ref{ex:23independenceb}, \ref{ex:toblerone} and \ref{ex:23independencec}.
Here ${\rm Hu}_A$ is a hyperdeterminant. So, we can use
the dedicated techniques described in \cite[Chapter 14]{GKZ}
to represent ${\rm Hu}_A$ and thereby compute~$\Lambda_A$.

A third technique for computing $\Lambda_A$ uses
the univariate polynomial $p_i(t_i)$ in (\ref{eq:principal}).
The logarithmic $A$-discriminant usually appears as
the largest factor in the discriminant of $p_i(t_i)$.

Our final method for computing $\Lambda_A$ rests on the
logarithmic Horn uniformization in Theorem \ref{thm:logpara}.
This gives a parametric representation of $\Lambda_A$, and
we can use specific techniques for the implicitization of
rational hypersurfaces. In our view, this fourth method is
the most important one. It is the foundation for
 the numerical methods we shall describe below.

To make Theorem \ref{thm:logpara} completely explicit, we
choose an $n\times(n-d-1)$ matrix $\hat B$ whose columns span the kernel of $\hat A$.
% such that $\hat A\hat B=0$.
We write  $\hat a_i \in \ZZ^{d+1}$ for the $i$th column of $\hat A$
and $\hat b_i$ for the $i$th row of $\hat B$. The free variables in our parametrization are
$u\in\CC^{d+1}$, $w\in\CC^{n-d-1}$, and $t\in(\CC^*)^d$. 
With these conventions, we write the generic point in $L^\perp \star L^{-1}$ as follows:
\begin{equation}
\label{eq:phi}
\phi\,\,=\,\,\left(\,
\frac{\hat b_1\cdot w}{\hat a_1\cdot u},\,
\frac{\hat b_2\cdot w}{\hat a_2\cdot u},\,
\dots\,,\,\frac{\hat b_n\cdot w}{\hat a_n\cdot u}\,\right).
\end{equation}
The resulting parametrization of $\Lambda_A$ by the unknowns
 $t,u,w$ equals
\begin{equation}
\label{eq:para1}
z\,=\, (t^{-a_1},\dots,t^{-a_n}) \star \phi  \quad {\rm and} \quad
s \,=\, \sum_{i=1}^n \phi_i \quad {\rm and} \quad  \nu\, =\, -A\phi.
\end{equation}

In the case $m \geq 2$, we have an analogous parametrization 
of $\Lambda_\boldA$, thanks to Theorem \ref{thm:hornm2}. We write $n=n_1+\cdots+n_m$.
Here $\hat A$ is replaced by the Cayley matrix ${\rm Cay}(\boldA)$, which has
 format $(m+d) \times n$, and $\hat B$ is 
an $n\times(n-d-m)$ matrix with linearly independent columns such that
 ${\rm Cay}(\boldA) \hat B = 0$.
 With this, we can now recycle the formulas  in (\ref{eq:phi}) and (\ref{eq:para1}).

\begin{example}[$m=d=2$] \label{ex:23independenced}
We revisit Examples~\ref{ex:twolines} 
and \ref{ex:23independencec}. Then
$\hat A = {\rm Cay}(\boldA) $ and
$$
\hat B^T \,\, = \,\, 
\begin{bmatrix}
                            \,\,  1 &  \! -1 &    0  &  -1  &  \,1\, &   0\,\, \\
                           \,\,   1  &   0  &  \!-1  &  -1  &  \,0\,  &  1\,\, \\ 
\end{bmatrix}.
$$    
Our parametrization of  the logarithmic $\boldA$-discriminant 
$\Lambda_\boldA$ uses the variables
$ (u_1,u_2,u_3,u_4)$, $ (w_1,w_2)$ and $ (y_1,y_2,t_1,t_2)$.
The general recipe in (\ref{eq:phi}) and (\ref{eq:para1}) translates into the formulas
\begin{small} $$
z_{11} \,=\, \,y_1 \cdot \frac{w_1+w_2}{u_1} \,,\,\quad
z_{12} \,=\, -\, y_1 t_1 \cdot \frac{w_1}{u_1+u_3}\,,\,\quad
z_{13}\, =\, - \,y_1 t_2 \cdot \frac{w_2}{u_1+u_4}\, ,
$$
$$
z_{21}\, =\, - \,y_2 \cdot \frac{w_1+w_2}{u_2}\, , \,\quad
z_{22} \, =\, \, y_2 t_1 \cdot  \frac{w_1}{u_2+u_3} \,, \,\quad
z_{23}\, =\, \,y_2 t_2 \cdot \frac{w_2}{u_2+u_4} \,, 
$$
$$
s_1 \,=\, \frac{w_1+w_2}{u_1} \,- \,\frac{w_1}{u_1+u_3} \,-\, \frac{w_2}{u_1+u_4} \,, \quad
s_2 \,=\, -\frac{w_1+w_2}{u_2}  \,+\, \frac{w_1}{u_2+u_3} \,+\, \frac{w_2}{u_2+u_4},
$$
$$
\nu_1 \,\,=\,\, \frac{w_1}{u_1+u_3} \,-\, \frac{w_1}{u_2+u_3} \,, \quad
\nu_2 \,\,=\,\, \frac{w_2}{u_1+u_4} \,-\, \frac{w_2}{u_2+u_4}.
$$ \end{small}
Using implicitization, we obtain the sum
of $305$ monomials of degree $(6,6,4) $ in
Example \ref{ex:twolines}.
\end{example}

All of these symbolic elimination techniques are costly. It is
infeasible to compute the full expansion of
$\Lambda_A$ for larger matrices $A$.
However, we can still compute a numerical representation of
 $\Lambda_A$ which yields the bidegree.
 The same works to compute the multidegree of $\Lambda_{\boldA}$ when $m\geq 2$ as in Section~\ref{sec7}.
We developed such an algorithm along with an implementation in the
software \texttt{HomotopyContinuation.jl} \cite{BT}.
In what follows, we explain our method.

We are given the generators of the ideal $I$ in (\ref{eq:idealtoelim}),
and we know that the elimination ideal
$I\cap\QQ[z,s,\nu] $ is principal.
However, we do not know its generator $\Lambda_A$.
Let $V(I)$ denote the incidence correspondence in the triple product
$(\CC^*)^{d+2} \times\PP^{n-1}\times\PP^d$ with 
coordinates $(t_1,\dots,t_d,y_1,y_2)$, $(z_1:\dots:z_n)$, and $(s:\nu_1:\dots:\nu_d)$. 
We wish to study the bidegree of the closure of the projection of $V(I)$ into $\PP^{n-1}\times\PP^d$, which is the logarithmic $A$-discriminant $\Lambda_A$.

For the symbolic elimination implemented in our \texttt{Oscar} snippet, 
the generators of $I$ must be polynomials. But, the software \texttt{HomotopyContinuation.jl} is 
capable of working with rational functions. So, we now use the
following  incidence variety  which is birational to $V(I)$:
$$\mathcal I\,=\,\{(t,z,(s,\nu))\in (\CC^*)^d\times\PP^{n-1}\times\PP^d\mid (t,z)\in\mathcal R, (s,-\nu)=\hat A_z\Phi_A(t)\}.$$

To compute the degree $d_z$ of $\Lambda_A$ in  $z$,
 we substitute $(s,\nu)$ in the equation of $\mathcal I$ by generic values $(s^*,\nu^*) \in \mathbb{P}^{d}$.
 We are left with polynomials in $t$ and $z$ which define an incidence variety ${\cal I}_{s^*, \nu^*} \subset (\CC^*)^{d} \times\PP^{n-1}$. The coordinate projection in $\PP^{n-1}$ is the hypersurface $\Lambda_A(z, s^*, \nu^*)$ obtained by specializing $\Lambda_A$, so it has degree $d_z$. Moreover, we know from Theorem \ref{thm:logpara} that, for a generic point $z^* \in \Lambda_A(z, s^*, \nu^*)$, there is precisely one degenerate critical point, so the coordinate projection to $\PP^{n-1}$ is generically one-to-one. We intersect ${\cal I}_{s^*, \nu^*}$ with $(\CC^*)^{d} \times L$, where $L \subseteq \mathbb{P}^{n-1}$ is a generic line. The number of intersection points is $d_z$. Computing the degree in $(s,\nu)$ is similar, by fixing $z = z^*$
and intersecting with a generic line in $(s, \nu)$-space.

In practice, we  compute intersections numerically  with \texttt{HomotopyContinuation.jl}. More precisely, we use the command \texttt{monodromy\_solve} to obtain all intersection points of ${\cal I}_{s^*, \nu^*} \cap ((\mathbb{C}^*)^{d} \times L)$.  We briefly explain how this works.  The coordinate projection  
\[ \pi\, : \, \mathcal{E}=\{ (t, z, (s,\nu), L) \in \mathcal I \times {\rm Gr}(2,n) \, : \, z \in L \} \, \rightarrow \, \PP^d \times {\rm Gr}(2, n) \]
is a branched covering of degree $d_z$. Suppose we are given a point $p^* = (t^*, z^*, (s^*,\nu^*), L^*)$ on  $\mathcal{E}$
and a closed loop $\gamma: [0, 1] \rightarrow \PP^d\times {\rm Gr}(2,n)$ such that $\gamma(0) = \gamma(1) = ((s^*,\nu^*),L^*)$. That loop lifts to a smooth \emph{solution path} on $\mathcal{E}$ which starts at $p^*$ and ends at a potentially different point in $\pi^{-1}((s^*,\nu^*),L^*)$. In this way, each such a loop induces a permutation of the $d_z$ points in the fiber $\pi^{-1}((s^*,\nu^*),L^*)$, and these permutations generate the monodromy group. Tracking solution paths amounts to solving an ODE, which is the method of \emph{homotopy continuation}. Repeating this process for sufficiently many random loops $\gamma$ is the \emph{monodromy method} \cite{Duff} for computing the full fiber $\pi^{-1}((s^*,\nu^*),L^*)$. This is very effective in~practice.

To start the monodromy method, we must compute the starting point $(t^*, z^*, (s^*,\nu^*), L^*)$, but this is trivial thanks to the logarithmic Horn uniformization. Indeed, Theorem \ref{thm:logpara} allows us to compute a generic point $(t^*,z^*,(s^*,\nu^*))$ on ${\cal I}$, and $L^*$ is a generic line passing through~$z^*$. 

When implementing this method, we work in affine charts of $\mathbb{P}^{n-1}, \mathbb{P}^d$ and ${\rm Gr}(2,n)$. For easier generalization to the $m>1$ case, we replace the matrix $\hat{A}_{z \star \Phi_A(t)} \hat{A}^T$ by the Jacobian matrix ${\rm Jac}(t)$ of the critical equations \eqref{eq:criteqs}. For computational efficiency, instead of requiring $\det({\rm Jac}(t))=0$, we require that ${\rm Jac}(t)$ has a nonzero kernel vector $q$. This introduces new indeterminates $q \in \mathbb{P}^d$, again in affine coordinates, 
which avoids expanding the determinant. For more implementation details, see the code on our Zenodo page \cite{Zenodo}.
This repository can be used to reproduce the computational results reported below. 

The numerical method  above works for
logarithmic $A$-discriminants whose degrees are in the thousands.
See Tables \ref{tab:logdet}, \ref{tab:cube}, \ref{tab:hypersimplex}
and \ref{tab:veronese} for some of those degrees.
If the toric variety $X_A$ is smooth, then we can try to
 compute the Hurwitz degree $h$ using the
Adjunction Formula. See
Corollary \ref{cor:logdet} and Proposition \ref{prop:veronese}.
By Theorem~\ref{thm:fromHurwitz}, we can then 
conclude that the bidegree is $(dh,h)$. In such cases,
theory and numerical computation validate each other.

\begin{table}[h!]
\centering
\begin{tabular}{|c|c|}
 \hline
 $d$ & bidegree \\ [0.5ex]
 \hline\hline
 $2$ & $ (4, 2)$  \\
 $3$ &$ (36, 12)$  \\
 $4$ &$ (288, 72)$  \\
 $5$ & $(2400, 480)$  \\ 
 $6$ & $(21600, 3600)$  \\ [1ex] 
  \hline
\end{tabular}
\caption{Bidegrees of logarithmic hyperdeterminants for $2 \!\times\! 2 \!\times\! \cdots \!\times\! 2$ tensors}
\label{tab:cube}
\end{table}

Table \ref{tab:cube} serves to illustrate this approach.
Here the support is the $d$-cube
$A = \{0,1\}^d$. Hence
 $n= 2^d$, $v=d\, !$, and $X_A$ is the
Segre embedding of $(\PP^1)^d$ into $\PP^{n-1}$.
In algebraic statistics, this is the binary independence model.
Using the Adjunction Formula, one can show that
$g=1+d!(d-3)/2$, and hence 
$h=2v+2g-2 = (d-1) \cdot d\,!$. 
The logarithmic $A$-discriminant $\Lambda_A$ for the cube $A$ was dubbed
{\em logarithmic hyperdeterminant} in Example \ref{ex:23independence}.

\begin{corollary} The logarithmic hyperdeterminant of a 
$2^d = 2 \!\times\! 2 \!\times \cdots \times\! 2$ tensor has bidegree
$$ (dh,h) \,\, = \,\,\bigl(\,d(d-1)\,d\,!\,,\, (d-1)\,d\,! \,\bigr). $$
\end{corollary}

Our numerical algorithm computes the bidegree of $\Lambda_A$ also when $X_A$ is not smooth.

\begin{example}\label{ex:hypersimplex}
Let $A$ be the {\em hypersimplex} $\Delta(r,s)$, i.e.~the convex hull of 
  all vectors  in $\{0,1\}^s$ that have coordinate sum $r$.
  This is a polytope of dimension $d=s-1$ with $\binom{s}{r}$ vertices.
   Table~\ref{tab:hypersimplex} shows the suspected bidegrees of $\Lambda_A$
   for some values of $r$ and $s$, computed by our numerical code at \cite{Zenodo}.
   Note that $b < dh$ holds for each table entry, so 
the inequality in
 Theorem~\ref{thm:fromHurwitz} is strict. This is consistent with the fact that $X_A$ is singular in all these cases. 
  \begin{table}[h!]
  \centering
  \begin{tabular}{|c||c|c|c|c|} 
  \hline
    $r \,\backslash \, s$ & 4 & 5 & 6 & 7 \\
  \hline\hline
    2 & $(12, 8)$ & $(80, 30)$ & $(330, 84)$   & $(1092, 210)$ \\
    3 &           & $(80, 30)$ & $(1080, 264)$ & $(7980, 1470)$ \\
    4 &           &            & $(330, 84)$   & $(7980, 1470)$ \\
    5 &           &            &               & $(1092, 210)$ \\
  \hline
  \end{tabular}
  \caption{Bidegree of the logarithmic $A$-discriminant when $A$ is the hypersimplex $\Delta(r,s)$.}
  \label{tab:hypersimplex}
  \end{table}
\end{example}
It remains a challenge to determine the bidegree of $\Lambda_A$ in cases where $X_A$ is singular. More precisely, what is the formula for $d_z=b$ in terms of the Hurwitz degree $h$, the degrees of the face discriminants $\Delta_{A_F}$, and the local multiplicities of $X_A$ along its torus orbits? Numerical computations like those in Example \ref{ex:hypersimplex} help to formulate and test conjectures.

\section{Statistics and Physics}
\label{sec9}
    
Our study of the logarithmic $A$-discriminant is motivated by
current research trends in algebraic statistics and particle physics.
It is our aim in this section to explain these applications.

We begin with algebraic statistics.
In that field, $\Lambda_A$ is a parametric version of the {\em data-discriminant}
due to Rodriguez and Tang \cite{RT}. The setting is the
   representation of scaled toric models  used in
\cite{ABB}, \cite[Section 16.1]{Telen} and \cite[Section 5]{TW}.
We fix a matrix $A = (a_1,a_2,\ldots,a_n) \in \NN^{d \times n}$ as above, and we consider
an arbitrary positive vector
$z \in \RR^n_{>0}$. The polynomial $f_z(t) = \sum_{i=1}^n z_i t^{a_i}$ 
takes on positive values for positive arguments $t \in \RR^d_{>0}$.

The toric model defined by the pair $(A,z)$ is the intersection of the
scaled toric variety
$z \star X_A$ with the $(n-1)$-dimensional probability simplex,
here identified with
$\PP^{n-1}_{> 0}$. This intersection is a statistical model
with state space $[n] = \{1,2,\ldots,n\}$. The probability of
observing state $i$ equals $z_i\, t^{a_i} / f_z(t)$.
Note that the $n$ probabilities are positive and they sum to $1$.
The Zariski closure of the model in $\PP^{n-1}$ is the
$d$-dimensional projective  variety $z \star X_A$.

Consider a data vector $u = (u_1,u_2,\ldots,u_n) \in \NN^n$.
The $i$th coordinate $u_i$ is the number of times 
the state $i$ was observed. The sum
$|u| = u_1+u_2 + \cdots + u_n$ is the sample size.
The likelihood function is the product of the probabilities of each of the outcomes 
in the sample:
\begin{equation}
\label{eq:loglike2} \binom{|u|}{u} \cdot \prod_{i=1}^n \biggl[\frac{z_i \, t^{a_i}}{f_z(t)} \biggr]^{u_i} \,\, = \,\,
z^u \, \binom{|u|}{u} \cdot f_z(t)^s t^\nu, \qquad \hbox{where $\,s = - |u|\,$ and $\,\nu = A u$.}
\end{equation}
The multinomial coefficient $\binom{|u|}{u}$ is a constant, as is $z^u$.
After deleting these constants, (\ref{eq:loglike2}) equals $f_z(t)^st^\nu$.
This is the likelihood function of the scaled toric model.
It is a function of the model parameters $t$, for a fixed
model $(A,z)$ and fixed data $u$. The likelihood 
depends on the data $u$ via the sufficient statistics, which are
our exponents $s = - |u|$ and $\nu = A u$.

The aim of likelihood inference is to determine the
model parameters $\hat t$ which maximize the function
(\ref{eq:loglike2}). One calls $\hat t$ the
{\em maximum likelihood estimator (MLE)}. It is an algebraic 
 function of both the data $u$ and the scaling parameters $z$,
 where now only $A$ is kept fixed. The discriminant of that
 algebraic function is precisely the logarithmic $A$-discriminant $\Lambda_A$.
 This can now be viewed 
  as a bihomogeneous polynomial in the $n+n$ unknowns $z$ and $u$. 
 
 In its statistical interpretation,
  Example \ref{ex:square} is the scaled version of the independence model for two binary random variables.
The data for this model is a $2 \times 2$ table of nonnegative integers
$u = \left [\begin{smallmatrix} u_1  & u_2 \\ u_3 & u_4 \end{smallmatrix} \right ]$.
The sufficient statistics are the row sums and column sums of the
$2 \times 2$ table $u$.
These four numbers are dependent: any three of them determines the fourth.
In the notation of the present paper, the sufficient statistics are given by the three exponents
\begin{equation}
\label{eq:suffstat}
\nu_1 = u_2+u_4\,, \,\,\nu_2 = u_3+u_4 \,\,\, {\rm and} \,\,\, s = -u_1-u_2-u_3-u_4.
\end{equation}
After this substitution, the logarithmic $A$-discriminant
$\Lambda_A$ is a polynomial with $12$ terms which
is quadratic in the data $u$ and quartic in the scaling parameters $z$.
We find that
$$ \Lambda_A \,=\,
(u_1-u_4)^2 z_2^2 z_3^2 \, + \, (u_2-u_3)^2 z_1^2 z_4^2\,+\,
2 (u_1 u_2 + u_1 u_3 + 2 u_1 u_4 + 2 u_2 u_3 + u_2 u_4 + u_3 u_4) z_1 z_2 z_3 z_4.
$$
This discriminant never vanishes when both the  parameters $z$ and
the data $u$ are positive.

Recall (e.g.~from \cite{ABB}) that the
ML degree drops from $2$ to $1$ when $z$ lies on the
principal $A$-determinant $E_A$. This special locus in $\PP^3$ includes the
surface defined by the $A$-discriminant
$\Delta_A = z_1 z_4 - z_2 z_3$. On that surface, the
logarithmic $A$-discriminant is a perfect square:
$$ \Lambda_A \,\,=\,\, \bigl(z_1 z_4 (u_1+u_2+u_3+u_4) \bigr)^2. $$

\begin{figure}[h]
\centering
 \includegraphics[width = 0.65 \linewidth]{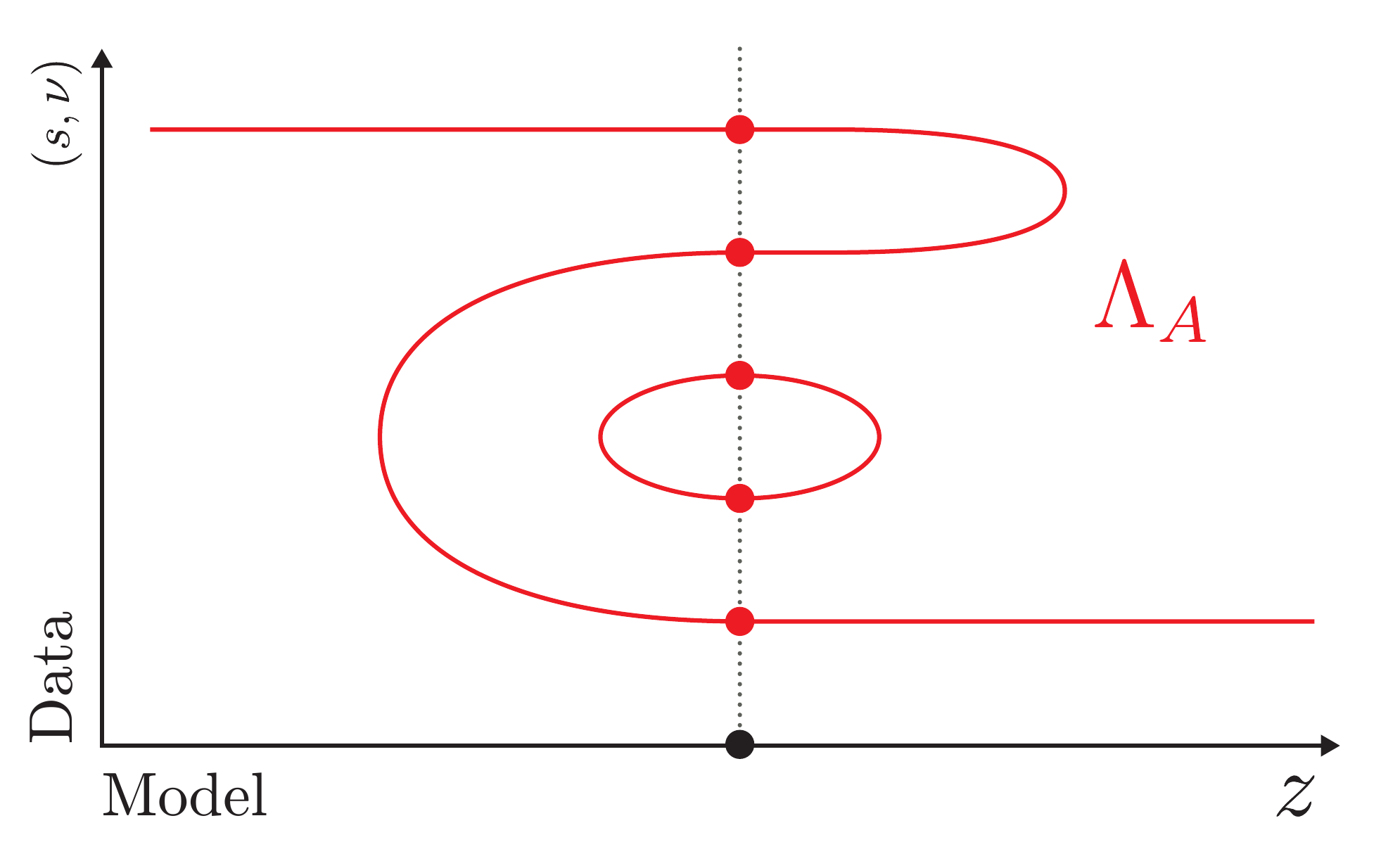} 
\caption{The logarithmic $A$-discriminant governs the relationship
between a toric model and its data. The former is fixed by a point on the
$z$-axis. The points on the curve over that model point $z$
reveal all sufficient statistics
$(s,\nu)$ which are degenerate for the chosen model.}
\label{fig:endler}
\end{figure}

Returning to the general  discussion, we have seen how $\Lambda_A$ extends the data-discriminant in \cite{RT}
by allowing the scaling parameters $z$ to vary. In our view, $\Lambda_A$ is a key tool for
studying degeneracies in the data locus of a family of toric models.
It is designed for parametric analysis  where  $z$ varies.
This variation can be visualized by a diagram as shown in
Figure~\ref{fig:endler}.

A staple of algebraic statistics is the {\em independence model}
for two random variables with $d_1$ and $d_2$ states respectively.
Geometrically, this is the Segre variety $X_A \simeq \PP^{d_1-1} \times
\PP^{d_2-1}$ in $\PP^{d_1 d_2-1}$.
The data is a  contingency table $u \in \NN^{d_1 \times d_2}$.
We just discussed the case $d_1 = d_2 = 2$.
A thorough study of all scalings of independence models
was undertaken by Clarke {\em et al.}~\cite{CHKO}.
The data-discriminant of these models 
is our logarithmic determinant.
Hence Corollary~\ref{cor:logdet}
and Table~\ref{tab:logdet} show the
degrees of the data-discriminant in
terms of both model $z$ and  data~$u$.

In Table~\ref{tab:cube} we increase the number of random variables from $2$ to $d$.
This shows the complexity of the data-discriminant of the independence model
for $d$ binary random variables.

We next consider independence models for identically distributed data.  Here
$X_A$ is the Veronese variety, comprising symmetric tensors of rank one, so our integrand in
(\ref{eq:eulerint}) is a general polynomial  $f_z$ of degree $r$ in $d$ variables.
We now present the analog to  Corollary~\ref{cor:logdet}.
For $r=2$, this gives the
logarithmic determinant of a symmetric $(d+1) \times (d+1)$ matrix. 

\begin{proposition} \label{prop:veronese}
Let $X_A \simeq \PP^d$ be the Veronese variety
of degree $v = r^d$ in $\PP^{\binom{d+r}{r}-1}$.
Then $\Lambda_A$ has bidegree $(dh,h)$, where $\, h = 2v + 2g - 2\,$ and
 $ \, g = 1 +  r^{d-1} ((d-1)r - d-1)/2$.
 \end{proposition}

\begin{proof}
We use Theorem~\ref{thm:fromHurwitz}.
The sectional genus $g$ of
$X_A$ is the genus of a curve in $\PP^d$
defined by $d-1$ general equations of degree $r$.
This is computed using the Adjunction Formula.
\end{proof}

Proposition \ref{prop:veronese} explains all entries in Table \ref{tab:veronese}.
The first row concerns univariate polynomials, as discussed in Example \ref{eq:rnc}.
 The entry for $d=2$ and $r=2$
appears in Example~\ref{ex:tact1}.

\begin{table}[h!]
\centering
\begin{tabular}{|c||c|c|c|c|c|c|}
 \hline
  $ d\,\backslash \, r $ & 2 & 3 & 4 & 5 & 6 & 7 \\
 \hline\hline
  $1$    & $(2,2)$   & $(4,4)$    & $(6,6)$    & $(8,8)$    & $(10,10)$  & $(12,12)$ \\
  $2$    & $(12,6)$  & $(36, 18)$ & $(72, 36)$ & $(120,60)$ & $(180,90)$ & $(252,126)$  \\
  $3$    & $(48,16) $&  $(216, 72)$ & $(576, 192)$ & $(1200, 400) $ & $(2160, 720)$ & $(3528, 1176)$ \\
  $4$    & (160,40) & $ (1080, 270) $ & $ (3840, 960) $ & $(10000, 2500) $ & $(21600, 5400)$ & $(41160, 10290)$ \\
   \hline
\end{tabular}
\caption{Bidegrees for the logarithmic $A$-discriminant
of a general polynomial of degree $r$}
\label{tab:veronese}
\end{table}

We conclude our discussion of algebraic statistics with an observation of general interest.

\begin{remark} \label{rmk:dependsonly} The expression for $\Lambda_A$ in $u$-coordinates 
 depends only on the toric variety $X_A$, whereas the $(s, \nu)$-representation depends on the specific choice of the integer matrix $A$. 
\end{remark}

We now shift gears, namely to the analysis of Feynman integrals in physics \cite{Weinzierl}.
One starts with a Feynman graph of genus $g$ with $d$ internal edges.
The Symanzik polynomials $\mathcal{U}$ and $\mathcal{F}$ of this graph are 
homogeneous of degree $g$ and $g+1$ respectively in the variables $t = (t_1,\ldots,t_d)$.
The Lee-Pomeransky representation of the Feynman integral takes the form (\ref{eq:eulerint}) 
where $f_z = \mathcal{F} + \mathcal{U}$. The coefficient vector 
is a linear function of physical parameters that represent masses and momenta of particles.
The support $A$ consists of certain monomials of degree $g$ and $g+1$, and 
one refers to $P = {\rm conv}(A)$ as the Feynman polytope of the graph.

If $g \geq 2$ then $z$ ranges over a low-dimensional subspace of $\RR^n$.
The principal $A$-determinant $E_A$ typically vanishes identically on that subspace, which indicates that a separate
theory of discriminants is needed to address this case. Such a  theory 
was developed 
 in \cite{FMT, MT}. For logarithmic discriminants, this is an 
interesting project left for future research.

The case of one-loop ($g=1$) is easier because here the GKZ setting
applies directly. This is explained in \cite[Section 4.1]{FMT}.
Here, $\mathcal{U} = \sum_{i=1}^d t_i$, and
$\mathcal{F}$ is a quadratic form in $t$ with general coefficients $z$.
The Newton polytope  is twice the unit $d$-simplex with a vertex cut~off:
$$ \qquad P \,=\,  {\rm conv} \bigl( e_1,e_2,\ldots,e_d,\,2e_1,2e_2,\ldots,2e_d \bigr) \qquad
{\rm and} \quad A \,=\, P \cap \NN^d. $$
The one-loop Feynman polytope $P$ has normalized volume $v = 2^d-1$,
and $n = \#(A) = d(d+3)/2$. The factorization of the principal 
$A$-determinant $E_A$ is given in \cite[Theorem~4.3]{FMT}.

\begin{example}[$d=2$] \label{ex:bubble} For the bubble diagram, we consider the $2 \times 5$ matrix $A=\left[\begin{smallmatrix}1 & 0 & 2 & 1 & 0 \\ 0 & 1 & 0 & 1 & 2\end{smallmatrix}\right]$. Our polynomial thus equals
$$ f_z(t) \,\, = \,\, z_1 t_1 + z_2 t_2 + z_3 t_1^2 + z_4 t_1 t_2 + z_5 t_2^2. $$
The principal $A$-determinant has degree $(d+1)v = 9$, and it
has six irreducible factors:
$$ E_A(z) \, = \, z_1 z_2 z_3 z_5 (4 z_3 z_5-z_4^2)(z_1^2 z_5-z_1 z_2 z_4+z_2^2 z_3) . $$
The last factor is the $A$-discriminant $\Delta_A$.
The logarithmic $A$-discriminant $\Lambda_A$ is a sum of $137$ monomials of bidegree~$(8,4)$.
Written as a quartic in $(s,\nu_1,\nu_2)$, it is dense with $15$ terms:
$$ \! \Lambda_A = 
(2 z_1 z_5-z_2 z_4)^2 (z_1 z_4-2 z_2 z_3)^2 s^4
\,+\,  4z_5 (z_1 z_4-z_2 z_3)^2 \Delta_A(z) \nu_1^4
\,+\,  4 z_3 (z_1 z_5-z_2 z_4)^2 \Delta_A(z) \nu_2^4 \,\,+\, \cdots
$$
The physics application \cite[Example 4.7]{FMT} uses a
Mandelstam variable $\sigma$ and masses $m_1,m_2$:
$$ z_1 = 1,\, \,z_2 = 1, \,\,z_3 = -m_1, \,\,z_4 = \sigma -m_1-m_2, \,\,z_5 = -m_2. $$
After this substitution, $\Lambda_A$ is a polynomial with $124$ monomials in
$ \sigma,m_1,m_2,s,\nu_1,\nu_2$.
\end{example}

For general dimensions $d$, we record the following corollary
to Theorem \ref{thm:fromHurwitz}.

\begin{corollary}
The logarithmic $A$-discriminant for the  Feynman polytope $P$ of the $d$-cycle
is a bihomogeneous polynomial $\Lambda_A$ of bidegree $(dh,h)$ in the 
unknowns $z$ and $(s,\nu)$, where
$$ h \,=\, 2^{d+1} + 2^{d-1}(d-3) - 2. $$
\end{corollary}

\begin{proof}
The variety $X_A$ is the blow-up of
$\PP^d$ at one point, so it is smooth.
Its sectional genus is the genus $g$ of 
the intersection of $d-1$ quadratic hypersurfaces in $\CC^d$. 
By setting $r=2$ in Proposition \ref{prop:veronese}, we find 
that $\, g \,=\, 1 + 2^{d-2} (d-3)$.
Now use the Hurwitz formula as before.
The degree can also be obtained by
subtracting $(2d,2)$ from	
 the $r=2$ column of
Table \ref{tab:veronese}.
\end{proof}

On first thought, one might believe that the logarithmic $A$-discriminant
plays a role in Landau analysis \cite{HPT, MT}. However, this is not the case. This
is explained by the factorization in Theorem \ref{thm:ERL}. 
Landau singularities are governed by the first factor $E_A(z)$.
The exponents $(s,\nu)$ are assumed to be non-resonant,
by analytic and dimensional regularization, so that $R_A(s,-\nu) \not= 0$.
The third factor $\Lambda_A(z,s,\nu)$ is not seen in Landau analysis. The number of \emph{master integrals}, i.e., the number of linearly independent integrals of the form \eqref{eq:eulerint} as the integration cycle $\Gamma$ varies, equals the number of 
critical points in $\mathcal{U}_z$ \emph{with multiplicities}.
If two critical points  of $f_z(t)^s t^\nu$ collide inside $\mathcal{U}_z$, 
then this does not create a Landau singularity.

The logarithmic $A$-discriminant is related to a more
refined wall crossing structure for integrals in physics.
This structure is discussed in recent work of 
Angius, Cacciatori and~Massidda \cite{ACM}.
Around each critical point of the integrand,
one defines a special integration contour $\Gamma$,
known as a {\em Lefschetz thimble}.
This works for generic parameters, and the
associated Feynman integrals \eqref{eq:eulerint} are linearly
independent. The Lefschetz thimbles span the twisted homology of $\,{\cal U}_z$ \cite[Section 3]{MMT}, and their integrals constitute a set of
master integrals. However, when two critical points
collide, this leads to a collapse in Lefschetz~thimbles.

Note that Lefschetz thimbles may fail to provide a basis for twisted homology even when the critical points do not collide. This is illustrated in \cite[Figure 6]{ACM}.
The correct \emph{analytic} discriminant locus consists of \emph{Stokes rays}, which 
are the walls to be crossed when the parameter $z$
varies in a one-dimensional complex family. 
This phenomenon is relevant for
 \emph{resurgence} in physics. The \emph{algebraic} condition $\Lambda_A = 0$ is called the \emph{vanishing discriminant case} in \cite[Figure 9]{ACM}. It is sufficient but not necessary for the Lefschetz thimbles to collapse.

Landau analysis focuses on objects
like $\Delta_A$ and $E_A$. These identify
potential singularities which can arise in the
$v$-dimensional family of Feynman integrals.
In contrast, our logarithmic $A$-discriminant $\Lambda_A$
offers more refined information about
the \emph{integrand}. It vanishes when $f_z(t)^s t^\nu$ is not a Morse function, which implies a
collapse in the thimble decomposition.

\bigskip \bigskip

\noindent {\bf Acknowledgement}:
We thank Sven Stawinski for pointing us to the reference \cite{ACM}.
 Our research was supported by
the European Research Council (ERC) through the synergy grant UNIVERSE+, 101118787.
$\!\!$ \begin{footnotesize}Views~and~opinions expressed
are however those of the authors only and do not necessarily reflect those of the European Union or the 
European
Research Council Executive Agency. Neither the European Union nor the granting authority
can be held responsible for them.
\end{footnotesize}

\bigskip \medskip

\noindent {\bf AI Declaration}:
LLMs helped us in deriving the inversion formula in Theorem \ref{thm:logpara}, in finding the deformation argument in Proposition \ref{prop:divides}, and in proofreading. We take full responsibility for the contents of this paper.

\bigskip
\bigskip \bigskip
		
\footnotesize
\noindent {\bf Authors' addresses:}
		
\smallskip
		
\noindent  Madeleine Goertz, UC Berkeley \hfill \url{m.goertz@berkeley.edu}
		
\noindent  Bernd Sturmfels, MPI-MiS Leipzig   \hfill \url{bernd@mis.mpg.de}

\noindent Simon Telen, MPI-MiS Leipzig \hfill \url{simon.telen@mis.mpg.de}

\end{document}